\documentclass[a4paper,11pt]{article}
\usepackage{amsmath}
\usepackage{amssymb}
\usepackage{amsthm}
\newcommand{\Hl}{{\mathbb H}}
\newcommand{\C}{{\mathbb C}}
\newcommand{\R}{{\mathbb R}}
\newcommand{\Q}{{\mathbb Q}}
\newcommand{\Z}{{\mathbb Z}}

\newcommand{\idP}{{\mathfrak P}}

\newcommand{\la}{{\mathfrak a}}

\newcommand{\mO}{{\cal O}}
\newcommand{\A}{{\mathbb A}}

\newcommand{\cross}{^{\times}}

\newcommand{\End}{\operatorname{End}}
\newcommand{\Tr}{\operatorname{Tr}}
\newcommand{\tr}{\operatorname{tr}}
\newcommand{\vol}{\operatorname{vol}}
\newcommand{\Ad}{\operatorname{Ad}}

\newcommand{\tk}{\tau_{\kappa}}
\newcommand{\ord}{\operatorname{ord}}

\newcommand{\ur}{\operatorname{ur}}
\newcommand{\spc}{\operatorname{sp}}
\newcommand{\fO}{{\mathfrak O}}
\newcommand{\cA}{{\cal A}}

\newcommand{\cL}{{\cal L}}
\newcommand{\cS}{{\cal S}}
\newcommand{\cT}{{\cal T}}

\newcommand{\JL}{\operatorname{JL}}

\newcommand{\diag}{\operatorname{diag}}
\numberwithin{equation}{section}
\theoremstyle{plain}
 \newtheorem{thm}{Theorem}[section]
 \newtheorem{prop}[thm]{Proposition}
 \newtheorem{lem}[thm]{Lemma}

 \newtheorem{assm}[thm]{Assumption}

\theoremstyle{definition}
 
 \newtheorem{rem}[thm]{Remark}

 \newtheorem{pf}{Proof}
 
\newcommand{\oo}{_{\infty}}

\begin{document}
\title{Fourier expansion of Arakawa lifting II:\\
Relation with central L-values}
\author{Atsushi Murase\thanks{Partly supported by Grand-in-Aid for Scientific Research (C) 20540031, Japan Society for the Promotion of Science and Kyoto Sangyo University Research Grants E1105.} \\
Department of Mathematical Science, Faculty of Science\\
Kyoto Sangyo University\\
Motoyama, Kamigamo, Kita-ku, Kyoto 603-8555, Japan\\ 
{\it E-mail address}:~murase@cc.kyoto-su.ac.jp\\[5pt]
Hiro-aki Narita\thanks{Partly supported by Grand-in-Aid for Young Scientists (B) 21740025, The Ministry of Education, Culture, Sports, Science and Technology, Japan.}\\
Graduate School of Science and Technology\\
Kumamoto University\\
Kurokami, Chuo-ku, Kumamoto 860-8555, Japan\\ 
{\it E-mail address}:~narita@sci.kumamoto-u.ac.jp}
\date{}
\maketitle
\begin{abstract}
This is a continuation of our previous paper [28]. The aim of the paper here is to study the Fourier coefficients of Arakawa lifts in relation with central values of automorphic $L$-functions. In the previous paper we provide an explicit formula for the Fourier coefficients in terms of toral integrals of automorphic forms on multiplicative groups of quaternion algebras. In this paper, after studying explicit relations between the toral integrals and the central $L$-values, we explicitly determine the constant of proportionality relating the square norm of a Fourier coefficient of an Arakawa lift with the central $L$-value. 
We can relate the square norm with the central value of some $L$-function of convolution type attached to the lift and a Hecke character. 
We also discuss the existence of strictly positive central values of the $L$-functions in our concern.\\[10pt]
{\it Key words}: Central $L$-values; Fourier coefficients; Quaternion unitary group; Theta lifts; Toral integrals.\\[10pt]
Mathematical Subject Classification 2010: 11F27,~11F30,~11F55,~11F67
\end{abstract}
\setcounter{section}{-1}
\section{Introduction.}
\subsection{}
Let $B$ be a definite quaternion algebra over $\Q$ with the discriminant $d_B$~(cf.~Section 1.1) and $D$ be a divisor of $d_B$. 
We recall that Arakawa lifting is a theta lifting to a cusp form on the quaternion unitary group $GSp(1,1)_{\A_{\Q}}$ from a pair consisting of an elliptic cusp form $f$ of level $D$ and an automorphic form $f'$ on $B^{\times}_{\A_{\Q}}$ ``with the same weight''~(for the definitions of the automorphic forms, see Section 1.1 and [28,~Sections 3.1 and 3.2]). 
At the archimedean place this theta lift $\cL(f,f')$ from $(f,f')$ generates a quaternionic discrete series representation in the sense of Gross-Wallach~[10]~(cf.~Section 1.1). 
Throughout the introduction we suppose that $(f,f')$ are non-zero Hecke eigenforms.
The Fourier expansion of a cusp form on $GSp(1,1)_{\A_{\Q}}$ is indexed by $\xi\in B^{-}\setminus\{0\}$ and a unitary character $\chi$ of $\A_{\Q}^{\times}\Q(\xi)^{\times}\backslash\A_{\Q(\xi)}^{\times}$~(cf.~Section 1.2), where $B^{-}$~(respectively $\Q(\xi)$) stands for the set of pure quaternions in $B$~(respectively the imaginary quadratic field generated by $\xi$ over $\Q$). 
We let $\cL(f,f')_{\xi}^{\chi}$ be the Fourier coefficient of $\cL(f,f')$ indexed by $\xi\in B^{-}\setminus\{0\}$ and $\chi$~(which is also called a Bessel period). 

As an application of our previous work [28] we study an explicit relation between the square norm of the Fourier coefficient $\cL(f,f')_{\xi}^{\chi}$ and the product of central values of the two quadratic base change lift $L$-functions for $(f,f')$ twisted by $\chi^{-1}$. 
We can furthermore relate the square norm with the central value of some automorphic $L$-function of convolution type attached to $\cL(f,f')$ and $\chi^{-1}$. 
Our method of the study here also yields an existence theorem of $\cL(f,f')$'s with strictly positive central values of the $L$-functions just mentioned. 
\subsection{}
Let $P_{\chi}(f;h)$ and $P_{\chi}(f';h')$ with $(h,h')\in GL_2(\A_{\Q})\times B^{\times}_{\A_{\Q}}$ be the toral integrals of $f$ and $f'$ defined by $\chi$ respectively~(cf.~Section 1.5). 
The main result [28,~Theorem 5.2.1] of our previous paper says that the Fourier coefficient $\cL(f,f')_{\xi}^{\chi}$ for a primitive $\xi\in B^{-}\setminus\{0\}$~(for the definition, see Section 1.3) is written as
\[
\cL(f,f')_{\xi}^{\chi}(g_{0})=C_0(f,f',\xi,\chi)\overline{P_{\chi}(f;\gamma_0)}P_{\chi}(f';\gamma'_0)
\]
with some $(g_0,\gamma_0,\gamma'_0)\in GSp(1,1)_{\A_{\Q}}\times GL_2(\A_{\Q})\times B_{\A_f}^{\times}$. 
Here the constant $C_0(f,f',\xi,\chi)$ depending on $f$, $f'$, $\chi$ and $\xi$ is explicitly determined~(see [28,~Theorem 5.2.1] or Theorem 1.3).

Let $\pi(f)$~(respectively~$\pi(f')$) be the automorphic representation generated by $f$~(respectively~$f'$) and let $\JL(\pi(f'))$ be the Jacquet-Langlands-Shimizu lift of $\pi(f')$~(cf.~[19],~[35,~Theorem 1]). 
We furthermore let $\Pi$~(respectively~$\Pi'$) be the base change lift of $\pi(f)$~(respectively~$\JL(\pi(f'))$) to $GL_2(\A_{\Q(\xi)})$. 
Denote the square norm of $f$~(respectively~$f'$) by $\langle f,f\rangle$~(respectively~$\langle f',f'\rangle$). 
Assuming that $f$ is primitive~(for the definition, see [25,~Section 4.6]), we now recall that, due to Waldspurger [38,~Proposition 7], there are constants $C(f,\chi)$ and $C(f',\chi)$ such that
\[
\frac{||P_{\chi}(f;\gamma_0)||^2}{\langle f,f\rangle}=C(f,\chi)L(\Pi,\chi^{-1},\frac{1}{2}),\quad\frac{||P_{\chi}(f';\gamma'_0)||^2}{\langle f',f'\rangle}=C(f',\chi)L(\Pi',\chi^{-1},\frac{1}{2}),
\]
where $L(\Pi,\chi^{-1},s)$~(respectively~$L(\Pi',\chi^{-1},s)$) denotes the $L$-function of $\Pi$~(respectively~$\Pi'$) with $\chi^{-1}$-twist. 
The results of Waldspurger [38,~Proposition 7] and ours~[28,~Theorem 5.2.1] thus imply that
\[
\frac{||\cL(f,f')_{\xi}^{\chi}(g_0)||^2}{\langle f,f\rangle\langle f',f'\rangle}=C(f,f',\xi,\chi)L(\Pi,\chi^{-1},\frac{1}{2})L(\Pi',\chi^{-1},\frac{1}{2}),
\]
where $C(f,f',\xi,\chi)$ is a constant depending only on $(f,f',\xi,\chi)$. 

We define some global spinor $L$-function for $GSp(1,1)_{\A_{\Q}}$ modifying Sugano's definition~(cf.~[36,~(3-4)]). 
More precisely, we define non-archimedean local factors of the $L$-function by the formula for the formal Hecke series and complete the global $L$-function with a suitable archimedean factor~(cf.~Section 2.6). 
This leads us to define the global $L$-function $L(\cL(f,f'),\chi^{-1},s)$ of convolution type for $\cL(f,f')$ and $\chi^{-1}$~(cf.~Section 2.6), whose local factors at unramified places are of degree eight. 
We see that this decomposes into
\[
L(\cL(f,f'),\chi^{-1},s)=L(\Pi,\chi^{-1},s)L(\Pi',\chi^{-1},s)
\]
(cf.~Proposition 2.10). 
In this paper we explicitly determine the constant $C(f,f',\xi,\chi)$~(cf.~Theorem 2.8) and deduce the following formula:
\begin{thm}[Theorem 2.11]
We have
\[
\frac{||\cL(f,f')_{\xi}^{\chi}(g_0)||^2}{\langle f,f\rangle\langle f',f'\rangle}=C(f,f',\xi,\chi)L(\cL(f,f'),\chi^{-1},\frac{1}{2}).
\]
\end{thm}
\subsection{}
Let us make several remarks on this theorem. 
We should first remark that Furusawa-Martin-Shalika [7]~[6] conjectured that the quantity $||\cL(f,f')_{\xi}^{\chi}(g_0)||^2$ is proportional to a similar $L$-function for the split symplectic group $GSp(2)$ of degree two with similitudes. 
Their conjecture is inspired by B{\"o}cherer's work~[2]. 
We note that the group $GSp(1,1)$ is an inner form of $GSp(2)$ and that the Langlands principle of functoriality~(cf.~[21]) suggests that an automorphic $L$-function for $GSp(1,1)$ should coincide with some automorphic $L$-function for $GSp(2)$. 

For any given divisor $D$ of $d_B$ we have a global maximal compact subgroup $K_f^D$~(cf.~Section 1.1). 
When $f$ is of level $D$, $\cL(f,f')$ is right $K_f^D$-invariant.
In the coming paper [30] the spinor $L$-function of $\cL(f,f')$ is proved to be that of  a paramodular cusp form on $GSp(2)_{\A_{\Q}}$ of level $d_BD$ given by some theta lift, where see Roberts and Schmidt [34,~Theorem 7.5.3,~Theorem 7.5.9,~Section A.6] for the non-archimedean local spinor $L$-functions of paramodular forms. 
For the case of $D=1$ this is essentially predicted by Ibukiyama~(cf.~[13], [14]). 
The $L$-function $L(\cL(f,f'),\chi^{-1},s)$ then turns out to be the $L$-function of convolution type attached to the above paramodular cusp form and $\chi^{-1}$.

For the determination of $C(f,f',\xi,\chi)$ we need explicit formulas for the two constants $C(f,\chi)$ and $C(f',\chi)$. 
There are many contributers to the study on $C(f,\chi)$ and $C(f',\chi)$, e.g. Gross [9], Hida [12], Martin-Whitehouse [24], Murase [26], Prasanna [32], Waldspurger [38], Xue [40] [41] and Zhang [42]. 
The works [9], [40] and [42] study the toral integrals under the geometric background. 
The method of [12], [26], [32], [38] and [41] is the theta correspondence, while that of [24] is the relative trace formula by Jacquet [17] and Jacquet-Chen [18,~Theorem 2]. 

We quote an explicit formula for $C(f,\chi)$ by [26]~(cf.~Proposition 2.6). 
To know $C(f',\chi)$ explicitly we recall that Waldspurger~[38,~Proposition 7] expresses it as a product of some local constants over places of $\Q$. More specifically, the local constants are written as products of some local integrals and  some ratios of local $L$-functions. 
To obtain an explicit form of $C(f',\chi)$ we have to evaluate the local integrals involved in $C(f',\chi)$. 
At a finite place not dividing the discriminant $d_B$ of $B$, we evaluate the integral by using Macdonald's explicit formula for zonal spherical functions~(cf.~[23,~Chap.V,~Section 3,~(3.4)]). The local integrals at other places are evaluated by a direct calculation. 
Our formula for $C(f',\chi)$ is stated as Proposition 2.7. 

We should remark that Martin and Whitehouse [24,~Theorem 4.1,~Theorem 4.2] already obtained similar formulas for $C(f,\chi)$ and $C(f',\chi)$. 
However, we note that 
the archimedean component of $f'$ has to be a highest weight vector~(which is not always true) in order to directly apply the formulas of [24] to $f'$. 
In addition, as far as we know, our method of the proof for Proposition 2.7 seems different from those of the known results.
\subsection{}
Let us specify the quaternion algebra as $B=\Q+\Q i+\Q j+\Q k$ with $i^2=j^2=-1$ and $ij=-ji=k$. 
As a primitive $\xi\in B^-\setminus\{0\}$ we take $\xi=i/2$. 
The level $D$ of $f$ is one or two since $d_B=2$ for this $B$. 
When $\chi$ is a unitary character of $\A_{\Q}^{\times}\Q(\xi)^{\times}\backslash\A_{\Q(\xi)}^{\times}$ unramified at every finite prime, we have verified the existence of simultaneously non-vanishing toral integrals for $(f,f')$~(cf.~[28,~Section 14] or Proposition 2.12). 
This and our explicit formulas for the toral integrals show the existence of simultaneously strictly positive central $L$-values
\[
L(\Pi,\chi^{-1},\frac{1}{2})>0,\quad L(\Pi',\chi^{-1},\frac{1}{2})>0
\]
(cf.~Theorem 2.13). 
For this result we remark that there are several results on the non-negativity of central values of these $L$-functions~(cf.~Guo [11],~Jacquet-Chen [18]). 
We furthermore see the existence of $\cL(f,f')$'s with the strictly positive central $L$-value as follows:
\begin{thm}[Proposition 2.12,~Theorem 2.14]
Let $B$ and $\chi$ be as above. When $D=1$~(respectively~$D=2$) 
let $\kappa\ge 12$~(respectively $8$) be an integer divisible by $4$~(respectively $8$). 
There exists $(f,f')$ of the same weight $\kappa$ such that 
\[
\cL(f,f')\not\equiv 0~\text{and}~L(\cL(f,f'),\chi^{-1},\frac{1}{2})>0.
\]
\end{thm}
We note that this $L$-function should be a Rankin-Selberg convolution $L$-function for $GSp(2)\times GL(2)$. 
As a related work we cite Lapid [22]. 
This deals with the non-negativity of central values of Rankin-Selberg L-functions for $SO(2m+1)\times GL(n)$, assuming that cuspidal automorphic representations of $GL(n)$~(respectively~$SO(2m+1)$) are ``orthogonal''~(respectively~generic). 
For this we note that $PGSp(2)\simeq SO(2,3)$.
\subsection{}
The outline of this paper is as follows. 
In Section 1 we review the main result of our previous paper [28]. In Section 2 we deal with our theorems: an explicit relation between $||\cL(f,f')_{\xi}^{\chi}(g_0)||^2$ and the central $L$-values, and an existence theorem of $\cL(f,f')$'s with the strictly positive central $L$-values. 
More precisely, after introducing basic notations of automorphic $L$-functions for $GL(2)$ in Sections 2.1 and 2.2, we give the explicit formula for $C(f,\chi)$ in Section 2.3 following [26]. 
In Section 2.4 we state our formula for $C(f',\chi)$. 
We then have our theorem relating $||\cL(f,f')_{\xi}^{\chi}(g_0)||^2$ to the product $L(\Pi,\chi^{-1},\frac{1}{2})L(\Pi',\chi^{-1},\frac{1}{2})$ of the central $L$-values in Section 2.5. 
In Section 2.6 we introduce the global $L$-function $L(\cL(f,f'),\chi^{-1},s)$ of convolution type attached to $\cL(f,f')$ and $\chi^{-1}$. 
In Theorem 2.11~(cf.~Section 2.6) we relate $||\cL(f,f')_{\xi}^{\chi}(g_0)||^2$ to $L(\cL(f,f'),\chi^{-1},\frac{1}{2})$. 
In Section 2.7 we show the existence of $\cL(f,f')$'s satisfying the strictly positivity of central values for the $L$-functions in our concern. 
We are then left with the proof of the formula for $C(f',\chi)$, which is considered in Section 3. 
After reviewing Waldspurger's formula [38,~Proposition 7] and Macdonald's formula [23,~Chap.V,~Section 3,~(3.4)] in Section 3.1 we evaluate each local component of $C(f',\chi)$. 
In Section 3.2~(respectively~Section 3.3) we deal with the local component at a finite prime $p\nmid d_B$~(respectively~$p|d_B$). 
In Section 3.4 we calculate the archimedean component of $C(f',\chi)$.
\subsection{}
(1)~In the introduction of our previous paper [28], we cited Jacquet [16] as a paper dealing with another proof of the formula of Waldspurger [37]~(not [38]). 
In fact, we can find in [16] an approach by a relative trace formula toward the formula. However, we remark that the proof of Waldspurger's formula [37] by the relative trace formula is not completed in [16]. 
We note that, for instance, Baruch-Mao [1] carry out the relative trace formula approach to the formula of Waldspurger [37].\\
(2)~In [31] we generalize our results without restricting ourselves to the case where the weights of $(f,f')$ are the same. 
However, we remark that this paper and our previous paper [28] include the essential part of the study at non-archimedean places necessary to obtain such generalized results in [31]. 
\subsection*{Notation.}
The ring $\A_{K}$ denotes the adele ring of a number field $K$ and $\A_f$ the ring of finite adeles in $\A_{\Q}$. 
For an algebraic group ${\cal G}$ over a field $F$ and an $F$-algebra $R$, ${\cal G}_R$ stands for the group of $R$-valued points. 
When $F=\Q$ and $R=\Q_p$ for a place $p$ of $\Q$ we sometimes denote ${\cal G}_{\Q_p}$ simply by ${\cal G}_p$. 
When $F$ is a number field or its completion at a finite place, $\mO_F$ denotes the ring of integers in $F$. 
For a finite set $S$, $|S|$ means the cardinality of $S$. 
Given a condition $C$, we put $\delta(C)=
\begin{cases}
1 & (\text{$C$ holds})\\
0 & (\text{otherwise})
\end{cases}$. 
For a measurable set $M$, $\vol(M)$ denotes its volume.
\section{Reviews on Arakawa lifting and its Fourier expansion.}
\subsection{}
In this section we basically use the notation in [28]. 
Let $B$ be a definite quaternion algebra over $\Q$. 
The discriminant $d_B$ of $B$ is defined as the product of finite prime $p$'s such that $B_p=B\otimes_{\Q}\Q_p$ is a division algebra. 
Let $B\ni x\mapsto \bar{x}\in B$ be the main involution of $B$ and denote by $n$ and $\tr$ the reduced norm and the reduced trace of $B$ respectively.
Namely
\[
n(x)=x\bar{x},\quad\tr(x)=x+\bar{x}
\]
for $x\in B$. 
For $X=
\begin{pmatrix}
a & b\\
c & d
\end{pmatrix}\in M_2(B)$ we put $\bar{X}=
\begin{pmatrix}
\bar{a} & \bar{b}\\
\bar{c} & \bar{d}
\end{pmatrix}$. 
Let $G=GSp(1,1)$ be the $\Q$-algebraic group defined by
\[
G_{\Q}:=\{g\in M_2(B)\mid {}^t\bar{g}Qg=\nu(g)Q,~\nu(g)\in\Q^{\times}\},
\]
where $Q:=
\begin{pmatrix}
0 & 1\\
1 & 0
\end{pmatrix}$. By $Z_G$ we denote the center of $G$.

We put $G_{\infty}^1:=\{g\in M_2(\Hl)\mid {}^t\bar{g}Qg=Q\}$, where $\Hl:=B{\otimes}_{\Q}\R$ is the Hamilton quaternion algebra. Then
\[
K_{\infty}:=\{
\begin{pmatrix}
a & b\\
b & a
\end{pmatrix}\in M_2(\Hl)\mid a\pm b\in\Hl^{1}\}
\]
forms a maximal compact subgroup of $G_{\infty}^1$, where $\Hl^1:=\{u\in\Hl\mid n(u)=1\}$. For a non-negative integer $\kappa$ we let $(\sigma'_{\kappa},V_{\kappa})$ be the $\kappa$-th symmetric tensor representation of $GL_2(\C)$ and $\sigma_{\kappa}$ the pull-back of $\sigma'_{\kappa}$ to $\Hl^{\times}$ via the standard embedding $\Hl\subset M_2(\C)$~(cf.~[28,~(1.4)]). 
This induces an irreducible representation $(\tk,V_{\kappa})$ of $K_{\infty}$ by
\[
\tk(
\begin{pmatrix}
a & b\\
b & a
\end{pmatrix}):=\sigma_{\kappa}(a+b)~(
\begin{pmatrix}
a & b\\
b & a
\end{pmatrix}\in K_{\infty}).
\]
Let $H$ and $H'$ be $\Q$-algebraic groups defined by
\[
H_{\Q}=GL_2(\Q),~H'_{\Q}:=B^{\times}
\]
respectively. 

Fix a maximal order $\fO$ of $B$ and a divisor $D$ of $d_B$. 
For $p|d_B$ let ${\idP}_p$ be the maximal ideal of the $p$-adic completion $\fO_p$ of $\fO$ and let
\[
L_p:=
\begin{cases}
{}^t(\fO_p\oplus\fO_p)&(\text{$p\nmid d_B$ or $p|D$}),\\
{}^t(\fO_p\oplus\idP_p^{-1})&(p|\frac{d_B}{D}).
\end{cases}
\]
We put $K_p:=\{k\in G_p\mid kL_p=L_p\}$ for each finite prime $p$ and $K_f^D:=\prod_{p<\infty}K_p$. 
The group $K_p$ forms a maximal compact subgroup of $G_p$ for each $p$ as is remarked in [27,~Section 2]. 

Hereafter $\kappa$ denotes an even integer. 
For $\kappa>4$ we then introduce the space $\cS_{\kappa}^D$ of $V_{\kappa}$-valued cusp forms $F$ on $G_{\A_{\Q}}$ satisfying the following:
\begin{enumerate}
\item $F(z\gamma gk_fk_{\infty})=\tk(k_{\infty})^{-1}F(g)$ for $(z,\gamma,g,k_f,k_{\infty})\in Z_{G,\A_{\Q}}\times G_{\Q}\times G_{\A_{\Q}}\times K_f^D\times K_{\infty}$,
\item for each fixed $g_f\in G_{\A_f}$, the right translations of the coefficients of $F|_{G_{\infty}^1}(g_f*)$ by elements in $G_{\infty}^1$ generate, as a $({\frak g},K_{\infty})$-module, the quaternionic discrete series representation~(cf.~[10]) with minimal $K_{\infty}$-type $\tk$~(if $F$ is non-zero). 
Here ${\frak g}$ denotes the Lie algebra of $G_{\infty}^1$.
\end{enumerate}
We remark that in [27,~Definition 2.1] and [28,~Definition 1.4.1] the second condition is replaced by the recurrence condition with respect to some reproducing kernel function. 
The equivalence of the two conditions is verified in [29,~Theorem 8.7,~Section 9].

For a positive integer $\kappa$ we let $S_{\kappa}(D)$ be the space of elliptic cusp forms of weight $\kappa$ with level $D$~(cf.~[28,~Section 3.1]) and ${\cal A}_{\kappa}$ be the space of automorphic forms of weight $\sigma_{\kappa}$ with respect to $\prod_{p<\infty}\fO_p^{\times}$~(cf.~[28,~Section 3.2]), where $\fO_p^{\times}$ denotes the unit group of $\fO_p$.
Now we can review the definition of Arakawa lifting. 
By using a metaplectic representation of $G_{\A_{\Q}}\times H_{\A_{\Q}}\times H'_{\A_{\Q}}$ we define in [27,~Section 4] the $\End(V_{\kappa})$-valued theta function $\theta^{\kappa}(g,h,h')$ with some specified $\End(V_{\kappa})$-valued Schwartz-Bruhat function on $B_{\A_{\Q}}^{\oplus 2}\times\A_{\Q}^{\times}$. 
Then, for $\kappa>4$, we have
\[
S_{\kappa}(D)\times {\cal A}_{\kappa}\ni (f,f')\mapsto {\cal L}(f,f')(g)\in\cS_{\kappa}^D
\]
(cf.~[28,~Theorem 3.3.2]) with
\[
{\cal L}(f,f')(g):=
\int_{(\R_+^{\times})^2(H\times H')_{\Q}\backslash (H\times H')_{\A_{\Q}}}\overline{f(h)}\theta^{\kappa}(g,h,h')f'(h')dhdh'.
\]
Here $(\R_+^{\times})^2$ means the connected component of the identity for the archimedean part of the center of $(H\times H')_{\A_{\Q}}$.
\subsection{}
We now review the Fourier expansion of ${\cal L}(f,f')$ described in [28,~Section 1.3]. We let $B^-:=\{x\in B\mid \tr(x)=0\}$ and have
\[
{\cal L}(f,f')(g)=\sum_{\xi\in B^-\setminus\{0\}}{\cal L}(f,f')_{\xi}(g),
\]
where
\[
{\cal L}(f,f')_{\xi}(g):=\int_{B^-\backslash B^-_{\A_{\Q}}}{\cal L}(f,f')(
\begin{pmatrix}
1 & x\\
0 & 1
\end{pmatrix}g)\psi(-\tr(\xi x))dx
\]
with the standard additive character $\psi:=\otimes_{p\le\infty}\psi_p$ of $\A_{\Q}/\Q$, which satisfies $\psi_{\infty}(x_{\infty})=\exp(2\pi\sqrt{-1}x_{\infty})$ for $x_{\infty}\in\R$. 
Here we normalize the measure $dx$ so that the volume of $B^-\backslash B^-_{\A_{\Q}}$ is one. For $\xi\in B^-\setminus\{0\}$ we let $X_{\xi}$ be the set of unitary characters on $\A_{\Q}^{\times}\Q(\xi)^{\times}\backslash \A_{\Q(\xi)}^{\times}$, which we call Hecke characters. 
With this $X_{\xi}$ the Fourier expansion is refined as follows:
\[
{\cL}(f,f')(g)=\sum_{\xi\in B^-\setminus\{0\}}\sum_{\chi\in X_{\xi}}{\cal L}(f,f')_{\xi}^{\chi}(g),
\]
where
\[
{\cal L}(f,f')_{\xi}^{\chi}(g):=\vol(\R_+^{\times}\A_{\Q}^{\times}\backslash\A_{\Q(\xi)}^{\times})^{-1}\int_{\R^{\times}_+\Q(\xi)^{\times}\backslash\A_{\Q(\xi)}^{\times}}{\cal L}(f,f')_{\xi}(s1_2\cdot g)\chi(s)^{-1}ds.
\]
\subsection{}
To review our explicit formula for ${\cL}(f,f')_{\xi}^{\chi}$ , we let $(f,f')\in S_{\kappa}(D)\times{\cal A}_{\kappa}$ and assume the following two conditions:\\
(1)~The two forms $f$ and $f'$ are Hecke eigenforms and have the same eigenvalue~(or signature) for the ``Atkin-Lehner involutions''. More precisely, for each $p|D$, let $\epsilon_p$~(respectively~$\epsilon'_p$) be the eigenvalue for the involutive action of 
$
\begin{pmatrix}
0 & 1\\
-p & 0
\end{pmatrix}
$~(respectively~a prime element $\varpi_{B,p}\in B_p$) on $f$~(respectively~$f'$). Then
\[
\epsilon_p=\epsilon'_p.
\]
Otherwise $\cL(f,f')\equiv 0$~(cf.~[27,~Remark 5.2 (ii)]).\\
(2)~We assume that $\xi\in B^-\setminus\{0\}$ is primitive~(cf.~[28,~Section 4.1]). Namely, for each finite prime $p$, we let
\[
\la_p:=
\begin{cases}
\fO_p&(\text{$p\nmid d_B$ or $p|D$})\\
\idP_p&(p|\frac{d_B}{D})
\end{cases},~
(\la_p^-)^*:=\{z\in B^-_{p}\mid\tr(\bar{z}w)\in\Z_p,~\text{for any $w\in\la_p\cap B_p^-$}\}
\]
and assume that
\[
\xi\in(\la_p^-)^*\setminus p(\la_p^-)^*.
\]
For this second assumption we note that, in general, the Fourier coefficient $F_{\xi}$ of an automorphic form $F$ on $G_{\A_{\Q}}$ satisfies
\[
F_{\xi}(g)=F_{\xi}(
\begin{pmatrix}
t & 0\\
0 & 1
\end{pmatrix}g)=F_{t\xi}(g)\quad(t\in\Q^{\times}).
\]
We then see that the problem determining $F_{\xi}$ is reduced to the case where $\xi$ is primitive.
\subsection{}
For each $\xi\in B^- \setminus \{0\}$, $d_{\xi}$ denotes the discriminant of $\Q(\xi)$. 
Hereafter, fixing $\xi$, we often denote $\Q(\xi)$ by $E$. 
We put
\[
a:=
\begin{cases}
2\sqrt{-n(\xi)}\sqrt{d_{\xi}}&(\text{$d_{\xi}$ is odd})\\
\sqrt{-n(\xi)}\sqrt{d_{\xi}}&(\text{$d_{\xi}$ is even})
\end{cases},\ b:=\xi^2-\frac{a^2}{4}.
\]
With these $a$ and $b$ we define $\iota_{\xi}:E^{\times}\hookrightarrow GL_2(\Q)$ by
\[
\iota_{\xi}(x+y\xi)=x\cdot 1_2+y\cdot
\begin{pmatrix}
\frac{a}{2} & b\\
1 & -\frac{a}{2}
\end{pmatrix}~(x,y\in\Q).
\]
Put $\theta:=r^{-1}(\xi-\frac{a}{2})$ with $r=\frac{2\sqrt{-n(\xi)}}{\sqrt{d_{\xi}}}\in\Q^{\times}$. Then $\{1,\theta\}$ forms a $\Z$-basis of the integer ring $\mO_E$ of $E$. We can rewrite $\iota_{\xi}$ as
\[
\iota_{\xi}(x+y\theta)=
\begin{pmatrix}
x & -rN_{E/\Q}(\theta)y\\
r^{-1}y & x+\Tr_{E/\Q}(\theta)y
\end{pmatrix}\quad(x,y\in\Q).
\]
The completion $E\oo$ of $E$ at $\infty$ is identified with $\C$ by
\[
\delta_{\xi}:E\oo\ni x+y\xi\mapsto x+y\sqrt{-n(\xi)}\in\C~(x,y\in \R).
\]
For a Hecke character $\chi=\prod_{p\leq\infty}\chi_p$ of $\R^{\times}_+E^{\times}\backslash\A_{E}^{\times}$, we let $w\oo(\chi)\in\Z$ be such that
\[
\chi\oo(u)=(\delta_{\xi}(u)/|\delta_{\xi}(u)|)^{w\oo(\chi)}~(u\in E\oo^{\times}).
\]
Furthermore, for each prime $p<\infty$, we let
\[
i_p(\chi):=\min\{i\ge 0\mid \chi|_{\mO_{E_p,i}^{\times}}\equiv 1\}
\]
with $\mO_{E_p,i}:=\Z_p+p^i\mO_{E_p}$ and let 
\[
\mu_p:=\frac{\ord_p(2\xi)^2-\ord_p(d_{\xi})}{2},
\]
which coincides with $\ord_p(r)$.
We now state the following~(cf.~[28,~Theorem 5.1.1]):
\begin{prop}
$\cL(f,f')_{\xi}^{\chi}\equiv 0$ unless $i_p(\chi)=0$ for any $p|d_B$ and $w_{\infty}(\chi)=-\kappa$.
\end{prop}
\subsection{}
In what follows, we assume the following condition on $\chi$:
\begin{assm}
$i_p(\chi)=0$ for any $p|d_B$ and $w_{\infty}(\chi)=-\kappa$.
\end{assm}
We need further notations to recall our formula for $\cL(f,f')_{\xi}^{\chi}$.
 We define $\gamma_0=(\gamma_{0,p})_{p\leq\infty}\in H_{\A_{\Q}}$ and $\gamma'_0=(\gamma'_{0,p})_{p<\infty}\in H'_{\A_f}$ as follows:
\begin{align*}
\gamma_{0,p}&:=
\begin{cases}
\begin{pmatrix}
1 & 0 \\
0 & p^{-\mu_p+i_p(\chi)}
\end{pmatrix}&(p\nmid D),\\
1_2 &\text{($p|D$ and $p$ is inert in $E$)},\\
\begin{pmatrix}
0 & 1\\
-p & 0
\end{pmatrix}&\text{($p|D$ and $p$ ramifies in $E$)},\\
\begin{pmatrix}
1 & a/2\\
0 & 1
\end{pmatrix}
\begin{pmatrix}
n(\xi)^{1/4} & 0\\
0 & n(\xi)^{-1/4}
\end{pmatrix}&(p=\infty),
\end{cases}\\
\gamma'_{0,p}&:=
\begin{cases}
\begin{pmatrix}
1 & 0\\
0 & p^{-\mu_p+i_p(\chi)}
\end{pmatrix}&(p\nmid d_B),\\
\varpi_{B,p}^{-1}&(p|d_B).
\end{cases}
\end{align*}
Here recall that $\varpi_{B,p}$ denotes a prime element of $B_p$~(cf.~Section 1.3).

In addition, we introduce the following local constants:
\[
C_p(f,\xi,\chi):=
\begin{cases}
p^{2\mu_p-i_p(\chi)}(1-\delta(i_p(\chi)>0)
e_p(E)p^{-1})&(p\nmid d_B),\\
1&(p|\frac{d_B}{D}),\\
2\epsilon_p&\text{($p|D$ and $p$ is inert in $E$)},\\
(p+1)^{-1}&\text{($p|D$ and $p$ ramifies in $E$)},
\end{cases}
\]
where
\[
 e_p(E)=
\begin{cases}
 -1&\text{($p$ is inert in $E$)},\\
 0&\text{($p$ ramifies in $E$)},\\
 1&\text{($p$ splits in $E$)}.
\end{cases}
\]
As in [28,~Section 2.4] we normalize the measure $ds=\prod_{p\le\infty}ds_p$ of $\A^{\times}_{E}$ so that
\[
\int_{\mO_{E_p}^{\times}}ds_p=1~\text{for any $p<\infty$},\quad\int_{E_{\infty}^1}ds_{\infty}=1,
\]
where $E_{\infty}^1:=\{s\in E_{\infty}^{\times}\mid n(s)=1\}$. 
In addition, we choose the measure of $\A_{\Q}^{\times}$ so that
\[
\vol(\Z_p^{\times})=1~\text{for any $p<\infty$}.
\]
For $(f,f')\in S_{\kappa}(D)\times\cA_{\kappa}$ we introduce their toral integrals
\[
P_{\chi}(f;h):=\int_{\R_+^{\times}E^{\times}\backslash \A_{E}^{\times}}f(\iota_{\xi}(s)h)\chi(s)^{-1}ds,~P_{\chi}(f';h'):=\int_{\R_+^{\times}E^{\times}\backslash \A_{E}^{\times}}f(sh')\chi(s)^{-1}ds
\]
with respect to $\chi$, 
where $(h,h')\in GL_2(\A_{\Q})\times B^{\times}_{\A_{\Q}}$. 
For these integrals we note that $E_p^{\times}$ is embedded into $B_p^{\times}\simeq GL_2(\Q_p)$ by $\iota_{\xi}$ when $p\nmid d_B$. 

We denote by ${\mathbf h}(E)$ and ${\mathbf w}(E)$ the class number of $E$ and the number of the roots of unity in $E$ respectively. Then we are able to state our formula for $\cL(f,f')_{\xi}^{\chi}$~(cf.~[28,~Theorem 5.2.1]).
\begin{thm}
Let $(f,f')$ be Hecke eigenforms with the same signature of the Atkin-Lehner involutions and $\xi\in B^-\setminus\{0\}$ be primitive~(cf.~Section 1.3 (1),~(2)), and recall Assumption 1.2 on $\chi$. 
We then have the following formula:
\begin{align*}
&\cL(f,f')_{\xi}^{\chi}(g_{0,f}
\begin{pmatrix}
\sqrt{\eta\oo} & 0\\
0 & \sqrt{\eta\oo}^{-1}
\end{pmatrix})\\
&=2^{\kappa-1}n(\xi)^{\kappa/4}\frac{\mathbf{w}(E)}{\mathbf{h}(E)}\cdot
\left(\prod_{p<\infty}C_p(f,\xi,\chi)\right)
\eta\oo^{\kappa/2+1}\exp(-4\pi\sqrt{n(\xi)}\eta\oo)
\,\overline{P_{\chi}(f;\gamma_0)}P_{\chi}(f';\gamma'_0).
\end{align*}
Here $\eta\oo\in \R\cross_+$ and $g_{0,f}=(g_{0,p})_{p<\infty}\in G_{\A_f}$ is given by
\[
g_{0,p}:=
\begin{cases}
\diag(p^{i_p(\chi)-\mu_p},p^{2(i_p(\chi)-\mu_p)},1,p^{i_p(\chi)-\mu_p})&(p\nmid d_B),\\
1_2&(p|d_B),
\end{cases}
\]
where the notation ``$\diag$'' means that $g_{0,p}$ is a diagonal matrix.
\end{thm}
\begin{rem}
According to Sugano [36,~Theorem 2-1], the Fourier coefficient ${\cal L}(f,f')_{\xi}^{\chi}$ is determined by the evaluation at $g_{0,f}
\begin{pmatrix}
\sqrt{\eta\oo} & 0\\
0 & \sqrt{\eta\oo}^{-1}
\end{pmatrix}$.
\end{rem}
\section{Relation with central $L$-values.}
\subsection{}
Throughout this section we let $(f,f')\in S_{\kappa}(D)\times {\cal A}_{\kappa}$ be Hecke eigenforms in the sense of [28,~Section 3.1,~3.2] and assume that $f$ is a primitive form~(for the definition, see [25,~Section 4.6]). 
Let $\pi(f)$~(respectively~$\pi(f')$) be the irreducible automorphic representation generated by $f$ (respectively~$f'$), and let $\JL(\pi(f'))$ be the Jacquet-Langlands-Shimizu lift of $\pi(f')$~(cf.~[19],~[35,~Theorem 1]). 
We remark that the classical prototype of the Jacquet-Langlands-Shimizu lift just mentioned is the Hecke equivariant isomorphism between ${\cal A}_{\kappa}$ and the space spanned by primitive forms in $S_{\kappa+2}(d_B)$, which is due to Eichler~(cf.~[4],~[5],~[35,~Section 6]). 
We denote by $\JL(f')$ the primitive form corresponding to $f'$ via this isomorphism. 

It is known that $\pi(f)$ and $\JL(\pi(f'))$ (respectively~$\pi(f')$) decompose into restricted tensor products over finite or infinite places $p\leq\infty$ of irreducible admissible representations of $GL_2(\Q_p)$~(respectively~$B_p^{\times}$). By $\pi_p$, $\pi'_p$ and $\pi''_p$ we denote the $p$-component of $\pi(f)$, $\JL(\pi(f'))$ and $\pi(f')$ respectively. According to such decompositions of $\pi(f)$ and $\pi(f')$, $f$ and $f'$ admit decompositions into pure tensor products
\[
\rho(f)=\prod_{p\leq\infty}f_p,~\rho'(f')=\prod_{p\leq\infty}f'_p,
\]
where we fix an isomorphism $\rho$~(respectively~$\rho'$) between $\pi(f)$ and the restricted tensor product $\otimes'_{p\le\infty}\pi_p$~(respectively~$\pi(f')$ and $\otimes'_{p\le\infty}\pi''_p$). 

We denote by $\Pi$~(respectively~$\Pi'$) the quadratic base change lift of $\pi(f)$~(respectively~$\JL(\pi(f'))$) to $GL_2(\A_{E})$. These $\Pi$ and $\Pi'$ also decompose into the restricted tensor products $\otimes'_{p\le\infty}\Pi_p$ and $\otimes'_{p\le\infty}\Pi'_p$ respectively, where each $\Pi_p$ or $\Pi'_p$ is a local base change lift of $\pi_p$ or $\pi'_p$ at every place $p$ respectively. We remark that each of the local and global representations just introduced has the trivial central character since so do $f$ and $f'$, thus it is self-dual~(cf.~[19,~Theorem 2.18 (i)],~[3,~Theorem 3.3.5]).
\subsection{Review on the adjoint $L$-functions and the $L$-functions of base change lifts for $GL_2$.}
Let $L(\pi,s)$ be the standard $L$-function for an automorphic representation $\pi$ of $GL_2(\A_{\Q})$. 
We denote by $L(\Pi,\chi^{-1},s)$~(respectively~$L(\Pi',\chi^{-1},s)$) the $L$-function of $\Pi$~(respectively~$\Pi'$) with $\chi^{-1}$-twist, and let $L(\pi(f),\Ad,s)$~(respectively~$L(\JL(\pi(f')),\Ad,s)$) be the adjoint $L$-function of $\pi(f)$~(respectively~$\JL(\pi(f'))$). 

We describe the local factors of $L(\Pi,\chi^{-1},s)$ and $L(\Pi',\chi^{-1},s)$~(respectively~$L(\pi(f),\Ad,s)$ and $L(\JL(\pi(f')),\Ad,s)$), following Jacquet~[15]~(respectively~Gelbart-Jacquet~[8]). We note that $\pi_p$~(respectively~$\pi'_p\simeq\pi''_p$) is a unitary unramified principal series representation of $GL_2(\Q_p)$ for each finite prime $p\nmid D$~(respectively~$p\nmid d_B$). 
This is due to the Ramanujan conjecture for holomorphic cusp forms on $GL_2$. 
In addition, we remark that, for $p|d_B$, $\pi''_p$ is written as
\[
B_p^{\times}\ni b\mapsto \delta_p\circ n(b)\in\{\pm 1\},
\]
with the unramified character $\delta_p$ of $\Q_p^{\times}$ of order at most two determined by $\delta_p(p)=\epsilon'_p$~(for $\epsilon'_p$ see Section 1.3). 
For $p|d_B$, $\pi'_p$ is thus the special representation of $GL_2(\Q_p)$ given by the irreducible subrepresentation~(or irreducible subquotient) 
of the normalized induced representation induced from the character of the Borel subgroup defined by two quasi-characters $\delta_p\cdot|*|_p^{\frac{1}{2}}$ and $\delta_p\cdot|*|_p^{-\frac{1}{2}}$~(cf.~[19,~Theorem 4.2 (iii)]), where $|*|_p$ is the $p$-adic valuation. 
We further note that $\pi_p$ is also a special representation when $p|D$. 
The induced representation giving rise to $\pi_p$ is associated with the unramified character $\delta'_p$ of order at most two determined by $\delta'_p(p)=-\epsilon_p$ in place of $\delta_p$, since the signature $\epsilon_p(=\epsilon'_p)$ of the Atkin Lehner involution of $f$ at $p|D$ differs from that of the corresponding automorphic form on $B^{\times}_{\A}$~(cf.~[19,~Lemma 15.7]). 
We remark that, when $p$ is inert in $E$ or when $p$ is ramified in $E$ and $\chi_p$ is unramified, $\chi_p$ can be written as
\[
\chi_p=\omega_p\circ n_{E_p/\Q_p}
\]
with a character $\omega_p$ of $\Q_p^{\times}$ of order at most two and the norm $n_{E_p/\Q_p}$ of $E_p$. In fact, $\omega_p\equiv1$ when $p$ is inert. 
In addition, at the archimedean place, $\pi_{\infty}$ and $\pi'_{\infty}$ are the discrete series representations with weight $\kappa$ and $\kappa+2$~(for the definition see [3,~Section 2.5]) respectively. 
For this fact see~[35,~Section 6].

We let $\pi_{\ur}$ be a unitary unramified principal series representation of $GL_2(\Q_p)$ with Satake parameter $(\alpha_p,\alpha_p^{-1})$ and the trivial central character, and let $\pi_{\spc}(+)$~(respectively $\pi_{\spc}(-)$) be the special representation $\pi'_p$ for $p|d_B$~(respectively $\pi_p$ for $p|D$). We denote by $\Pi_{\ur}$~(respectively~$\Pi_{\spc}(\pm)$) the base change lift of $\pi_{\ur}$~(respectively~$\pi_{\spc}(\pm)$) to $GL_2(E_p)$.

We first deal with the local $L$-functions of $\pi_{\ur}$ and $\Pi_{\ur}$. The following lemma is well-known.
\begin{lem}
(1)
\[
L_p(\pi_{\ur},s)=(1-\alpha_pp^{-s})^{-1}(1-\alpha_p^{-1}p^{-s})^{-1}.
\]
(2)
\[
L_p(\pi_{\ur},\Ad,s)=(1-p^{-s})^{-1}(1-\alpha_p^2p^{-s})^{-1}(1-\alpha_p^{-2}p^{-s})^{-1}.
\]
(3)
\begin{align*}
&L_p(\Pi_{\ur},\chi_p^{-1},s)=\\
&\begin{cases}
\prod_{i=1,2}(1-\alpha_p\chi_p(\varpi_{p,i})^{-1}p^{-s})^{-1}(1-\alpha_p^{-1}\chi_p(\varpi_{p,i})^{-1}p^{-s})^{-1}&(\text{$p$:split},~i_p(\chi)=0),\\
(1-\alpha_p^2p^{-2s})^{-1}(1-\alpha_p^{-2}p^{-2s})^{-1}&(\text{$p$:inert},~i_p(\chi)=0),\\
(1-\alpha_p\chi_p(\varpi_p)^{-1}p^{-s})^{-1}(1-\alpha_p^{-1}\chi_p(\varpi_p)^{-1}p^{-s})^{-1}&(\text{$p$:ramified},~i_p(\chi)=0),\\
1&(i_p(\chi)>0),
\end{cases}
\end{align*}
where $\varpi_{p,i}\in E_p$ with $i=1,2$~(respectively~$\varpi_p\in E_p$) denote prime elements of two distinct prime ideals dividing $p$~(respectively~a prime element dividing $p$) when $p$ is split~(respectively~$p$ is ramified).
\end{lem}
We next deal with the case of $\pi_{\spc}(\pm)$ and $\Pi_{\spc}(\pm)$, where note that $\pi_{\spc}(+)$ and $\Pi_{\spc}(+)$~(respectively~$\pi_{\spc}(-)$ and $\Pi_{\spc}(-)$) are defined for $p|d_B$~(respectively~$p|D$). 
\begin{lem}
We have
\begin{align*}
L_p(\pi_{\spc}(\pm),s)&=(1\mp\delta_p(p)p^{-(s+\frac{1}{2})})^{-1},\\
L_p(\pi_{\spc}(\pm),\Ad,s)&=(1-p^{-(s+1)})^{-1},\\
L_p(\Pi_{\spc}(\pm),\chi_p^{-1},s)&=
\begin{cases}
(1-p^{-(2s+1)})^{-1}&(\text{$p$:inert}),\\
(1\mp\delta_p(p)\omega_p(p)p^{-(s+\frac{1}{2})})^{-1}&(\text{$p$:ramified}).
\end{cases}
\end{align*}
\end{lem}
\begin{pf}
For the first formula see [19,~Proposition 3.6]. 
To verify the other two we need the local Rankin-Selberg convolution $L$-function $L_p(\pi_1\times\pi_2,s)$ of two irreducible admissible representations $\pi_1$ and $\pi_2$ of $GL_2(\Q_p)$~(cf.~[15]).

According to [8,~Proposition (1.4), (1.4.3)], we have
\[
L_p(\sigma_1\times\sigma_2,s)=L_p(\mu_1\mu_2,s)L_p(\nu_1\mu_2,s)
\]
for two special representations $\sigma_1$ and $\sigma_2$, 
where, for $i=1$ or $2$, $\sigma_i$ is attached to two quasi-character $\mu_i$ and $\nu_i$ of $\Q_p^{\times}$ with $\frac{\mu_i}{\nu_i}=|*|_p$. 
Then the second formula follows from this and the definition of the adjoint $L$-function in [8,~p.485], where we note that $\pi_{\spc}(\pm)$ is self-dual~(cf.~Section 2.1).

From [15,~Definition 20.1] we recall that
\[
L_p(\Pi_{\spc}(\pm),\chi_p^{-1},s)=L_p(\pi_{\spc}(\pm)\times\pi(\chi_p^{-1}),s),
\]
where see [19,~Theorem 4.6] for $\pi(\chi_p^{-1})$~(the dihedral representation associated with $\chi_p^{-1}$). 
From [15,~(20.3)] we can deduce the last formula.\qed
\end{pf}
For a positive integer $k\ge 2$ we let $\pi_{k}$ be the discrete series representation of $GL_2(\R)$ with weight $k$~(cf.~[3,~Section 2.5]) and $\Pi_k$ denote its base change lift. 
For $l\in\frac{1}{2}\Z$ we put $\chi_l$ to be the character of $\C^{\times}$ given by
\[
\chi_l(z):=\left(\frac{z}{\bar{z}}\right)^l,\quad z\in\C^{\times}.
\]
With 
\[
\Gamma_{\R}(s):=\pi^{-\frac{s}{2}}\Gamma(\frac{s}{2}),\quad\Gamma_{\C}(s):=2(2\pi)^{-s}\Gamma(s)
\]
we state the following:
\begin{lem}
We have
\begin{align*}
L_{\infty}(\pi_k,s)&=\Gamma_{\C}(s+\frac{k-1}{2}),\\
L_{\infty}(\pi_k,\Ad,s)&=\Gamma_{\R}(s+1)\Gamma_{\C}(s+k-1),\\
L_{\infty}(\Pi_k,\chi_l,s)&=
\begin{cases}
\Gamma_{\C}(s+\frac{k-1}{2}+|l|)\Gamma_{\C}(s-\frac{k-1}{2}+|l|)&(|l|\geq\frac{k-1}{2})\\
\Gamma_{\C}(s+\frac{k-1}{2}+|l|)\Gamma_{\C}(s+\frac{k-1}{2}-|l|)&(|l|\leq\frac{k-1}{2})
\end{cases},
\end{align*}
where $l\in\frac{1}{2}\Z$.
\end{lem}
\begin{pf}
We recall that the discrete series representation $\pi_k$, which is called a special representation in [19], corresponds to the representation of the Weil group $W_{\R}$ of $\R$ induced from the character $\chi_{\frac{k-1}{2}}$ of $\C^{\times}$~(cf.~[19,~Sections 5 and 12]), where note that $\C^{\times}$ is the Weil group of $\C$~(for the definition of Weil groups, see [39] and [19,~Section 12]).

For the first formula see [19,~p.181,~p.194,~p.195]. 
For the other two we need to consider the tensor product of two representations of $W_{\R}$ induced from characters of $\C^{\times}$, for which we refer to [15,~Case 17.6]. 
For the second formula recall the definition of the local adjoint L-function~(cf.~[8,~p.485]). 
We then see that
\[
L_{\infty}(\pi_{k},\Ad,s)=\Gamma_{\R}(s+1)\Gamma_{\C}(s+k-1).
\]
As for the last one, following the definition [15,~Definition 20.1] of the local L-factor of the lifting to $GL_2$ over a quadratic extension, we have $L_{\infty}(\Pi_{k},\chi_l,s)$ as above.\qed
\end{pf}
\subsection{Relation between $P_{\chi}(f;\gamma_0)$ and $L(\Pi,\chi^{-1},\frac{1}{2})$.}
By $\eta=\prod_{p\le\infty}\eta_p$ we denote the quadratic character attached to the quadratic extension $E/\Q$. We let $L(\eta,s)$ be the $L$-function defined by $\eta$ and $L_p(\eta_p,s)$ the local factor of $L(\eta,s)$ at a place $p$. 
We introduce the subsets $S_1$ and $S_2^{\pm}(f,\chi)$ of $d(D):=\{p<\infty\mid p|D\}$ which are denoted by ``$S_1(\xi)$'' and ``$S_2^{\pm}(f,\xi)$'' in [26,~Section 1.2] respectively. 
In view of [26,~Section 3.3] they are given as
\begin{align*}
S_1&:=\{p<\infty\mid \text{$p|D$, $p$ is inert in $E$}\},\\
S_2^{\pm}(f,\chi)&:=\{p<\infty\mid \text{$p|D$, $p$ ramifies in $E$, $\chi_p(\varpi_p)=\mp \epsilon_p$}\},
\end{align*}
where recall that $\epsilon_p$ denotes the eigenvalue of the Atkin-Lehner involution for $f$~(cf.~Section 1.3) and that $\varpi_p$ is a prime element of $E_p$~(cf.~Lemma 2.1). 
Note that, different from [26,~Section 1.2], $i_p(\chi)=0$ is already assumed for $p|D$~(cf.~Assumption 1.2) and that $S_1\cup S_2^+(f,\chi)\cup S_2^-(f,\chi)$ coincides with $d(D)$. 
We furthermore put $A(\chi):=\prod_{p<\infty}p^{i_p(\chi)}$.

We first quote [26,~Theorem 1.1], with a modification, for our situation. 
To consider the toral integral in [26,~(1.6)] we need another embedding $\iota'_{\xi}:E\hookrightarrow GL_2(\Q)$:
\[
\iota'_{\xi}(x+y\theta')=
\begin{pmatrix}
x & n(\theta')y\\
-y & x+\Tr(\theta')y
\end{pmatrix}.
\]
Here we put $\theta'=-\bar{\theta}$, 
for which we have to see that $\theta'$ satisfies the condition in [26,~Section 1.2]~(for this see also [28,~Lemma 4.3.1 (iii) and (4.9)]). 
Our embedding $\iota_{\xi}$ is related to $\iota'_{\xi}$ by
\[
\iota_{\xi}(\overline{x+y\theta})=
\begin{pmatrix}
r & 0\\
0 & 1
\end{pmatrix}\iota'_{\xi}(x+y\theta)
\begin{pmatrix}
r^{-1} & 0\\
0 & 1
\end{pmatrix}\quad(x+y\theta\in E),
\]
where recall that $r=\frac{2\sqrt{-n(\xi)}}{\sqrt{d_{\xi}}}$~(cf.~Section 1.4). 
Furthermore we note that our normalized measure of $\A_{E}^{\times}$ is the multiple of that of [26,~Section 2.4] by $\displaystyle\frac{\sqrt{|d_{\xi}|}}{4\pi}$. 
In fact, the volume of $\R^{\times}_+E^{\times}\backslash\A_E^{\times}$ with respect to our measure is equal to $\displaystyle\frac{{\bf h}(E)}{{\bf w}(E)}$ while that with respect to the measure of [26,~Section 2.4] is $2\prod_{p<\infty}L_p(\eta_p,1)=\displaystyle\frac{4\pi{\bf h}(E)}{\sqrt{|d_{\xi}|}{\bf w}(E)}$, where we note the normalization of the measure of $\A_{\Q}^{\times}$~(cf.~Section 1.5) to calculate this volume. 

We can then reformulate [26,~Theorem 1.1] as follows:
\begin{prop}
Let $L^{(\infty)}(\Pi,\chi^{-1},s)$ be the partial $\chi^{-1}$-twisted L-function of $\Pi$ without the archimedean factor. 
Under Assumption 1.2 we obtain
\[
|P_{\chi}(f;\gamma_0)|^2=
\begin{cases}
C'(f,\chi) L^{(\infty)}(\Pi,\chi^{-1},\frac{1}{2})&(S_1=S_2^+(f,\chi)=\emptyset),\\
0&(\text{otherwise}),
\end{cases}
\]
where
\[
C'(f,\chi):=
|d_{\xi}|(4\pi)^{-1-\kappa}(\kappa-1)!D^{-\frac{1}{2}}A(\chi)^{-1}2^{|d(D)|}\prod_{p|A(\chi)}L_p(\eta_p,1)^2.
\]
\end{prop}
Let $Z$ denote the center of $GL_2$. 
With the invariant measure $dg$ of $GL_2(\A_{\Q})$ as in [26,~Section 2.4], 
we introduce the Petersson norm
\[
\langle f,f\rangle:=\int_{Z_{\A_{\Q}}GL_2(\Q)\backslash GL_2(\A_{\Q})}|f(g)|^2dg
\]
of $f$. 
Rankin [33,~Theorem 3 (iii)] relates $\langle f,f\rangle$ to $L(\pi(f),\Ad,1)$ by taking the residue of the Rankin-Selberg integral of $f\times\bar{f}$ against an Eisenstein series. 
We can modify Rankin's argument adelically in a standard way to have the following:
\begin{prop}
\[
\langle f,f\rangle=2^{-\kappa-1}D\cdot L(\pi(f),\Ad,1).
\]
\end{prop}
In addition, we note that the archimedean factor of $L(\Pi,\chi^{-1},\frac{1}{2})$ is equal to
\[
\Gamma_{\C}(\kappa)\Gamma_{\C}(1),
\]
with $\Gamma_{\C}(s):=2(2\pi)^{-s}\Gamma(s)$~(cf.~Lemma 2.3). 
Then we can restate Proposition 2.4 as follows:
\begin{prop}
Under Assumption 1.2, for a non-zero primitive form $f$, we have
\[
\frac{|P_{\chi}(f;\gamma_0)|^2}{\langle f,f\rangle}=
\begin{cases}
C(f,\chi)L(\Pi,\chi^{-1},\frac{1}{2})&(S_1=S_2^+(f,\chi)=\emptyset),\\
0&(\text{otherwise})
\end{cases}
\]
with
\[
C(f,\chi):=
\frac{2^{|d(D)|-2}|d_{\xi}|\prod_{p|A(\chi)}L_p(\eta_p,1)^2}{D^{\frac{3}{2}}A(\chi)L(\pi(f),\rm{Ad},1)}.
\]
\end{prop}
\subsection{Relation between $P_{\chi}(f';\gamma'_0)$ and $L(\Pi',\chi^{-1},\frac{1}{2})$.}
Let $(*,*)_{\kappa}$ be a unitary inner product of $(\sigma_{\kappa}|_{\Hl^1},V_{\kappa})$~(for $\sigma_{\kappa}$ see Section 1.1) and $||v||$ denote the norm of $v\in V_{\kappa}$ induced by this inner product. 
We next provide an explicit relation between $||P_{\chi}(f';\gamma'_0)||^2$ and the central L-value $L(\Pi',\chi^{-1},\frac{1}{2})$. 
We postpone the proof until Section 3. 
To write down the relation we need several notations. 
We denote by $r_p$ the ramification index of $p$ for the quadratic extension $E/\Q$, i.e. $r_p:=
\begin{cases}
1&(\text{$p$:split or inert})\\
2&(\text{$p$:ramified})
\end{cases}$. 

Let $Z'$ denote the center of $B^{\times}$. 
We normalize the invariant measure $db:=\prod_{p\le\infty}db_p$ of $B_{\A_{\Q}}^{\times}$ so that
\[
\int_{\mO_p^{\times}}db_p=1~(p<\infty),\quad\int_{\Hl^{1}}db_{\infty}=1.
\]
We set
\[
\langle f',f'\rangle:=\int_{Z'_{\A_{\Q}}B^{\times}\backslash B^{\times}_{\A_{\Q}}}(f'(b),f'(b))_{\kappa}db.
\] 

We now note that the archimedean component $\pi''_{\infty}$ of $\pi(f')$ can be identified with an irreducible representation $(\sigma_{\kappa}|_{\Hl^1/\{\pm 1\}},V_{\kappa})$ of $B_{\infty}^{\times}/Z'_{\infty}\simeq\Hl^1/\{\pm 1\}$, where note that $\kappa$ is even~(cf.~Section 1.1) and $\sigma_{\kappa}$ is thus trivial on $\{\pm 1\}$. 
Let $v_{\kappa,\xi}$ be a highest weight vector of $V_{\kappa}$ with respect to $\sigma_{\kappa}(\R(\xi)^{\times})$-action, and let $v_{\kappa,\xi}^*\in V_{\kappa}$ be the dual of $v_{\kappa,\xi}$ with respect to $(*,*)_{\kappa}$. 
We set $f'_{\infty,\kappa}:=(f'_{\infty},v_{\kappa,\xi}^*)_{\kappa}v_{\kappa,\xi}$, where note that $f'_{\infty}$~(cf.~Section 2.1) is viewed as an element in $V_{\kappa}$. 
Then we have the following:
\begin{prop}
Under Assumption 1.2, for a non-zero Hecke eigenform $f'$, we have
\[
\frac{||P_{\chi}(f';\gamma'_0)||^2}{\langle f',f'\rangle}=
\begin{cases}
C(f',\chi)L(\Pi',\chi^{-1},\frac{1}{2})&(\text{$\pi''_p|_{E_p^{\times}}=\chi_p$ when $p$ divides $d_B$ and is ramified in $E$}),\\
0&(\text{otherwise}),
\end{cases}
\]
where
\[
C(f',\chi):=\frac{\sqrt{|d_{\xi}|}(\kappa+1)}{4A(\chi)L(\JL(\pi(f')),\Ad,1)}\prod_{p|A(\chi)}L_p(\eta_p,1)^2\cdot\prod_{p|d_B}r_pp^{-1}\cdot\frac{(f'_{\infty,\kappa},f'_{\infty,\kappa})_{\kappa}}{(f'_{\infty},f'_{\infty})_{\kappa}}.
\]
\end{prop}
\subsection{Main result~(the first form).}
Theorem 1.3, Proposition 2.6 and Proposition 2.7 imply the following theorem:
\begin{thm}
Under the assumption in Theorem 1.3 we have
\[
\frac{||\cL(f,f')_{\xi}^{\chi}(g_0)||^2}{\langle f,f\rangle\langle f',f'\rangle}=C(f,f',\xi,\chi)L(\Pi,\chi^{-1},\frac{1}{2})L(\Pi',\chi^{-1},\frac{1}{2}),
\]
where, 
if $\pi''_p|_{E_p^{\times}}=\chi_p$ for $p|d_B$ ramified in $E$ and $S_1=S_2^+(f,\chi)=\emptyset$,
\begin{align*}
&C(f,f',\xi,\chi)=\frac{2^{2\kappa+|d(D)|-6}(\kappa+1) n(\xi)^{\frac{\kappa}{2}}|d_{\xi}|^{\frac{3}{2}}\mathbf{w}(E)^2}{\mathbf{h}(E)^2A(\chi)^4D^{\frac{3}{2}}L(\pi(f),\Ad,1)L(\JL(\pi(f')),\Ad,1)}\\
&\prod_{p\nmid d_B}p^{4\mu_p}(1-\delta(i_p(\chi)>0)e_p(E)p^{-1})^2\prod_{p|d_B}r_pp^{-1}\prod_{p|D}(p+1)^{-2}\cdot e^{-8\pi\sqrt{n(\xi)}}\cdot\frac{(f'_{\infty,\kappa},f'_{\infty,\kappa})_{\kappa}}{(f'_{\infty},f'_{\infty})_{\kappa}}, 
\end{align*}
and $C(f,f',\xi,\chi)=0$ otherwise.
\end{thm}
\subsection{Main result~(the second form).}
We introduce the global $L$-function of convolution type attached to $\cL(f,f')$ and $\chi^{-1}$, and relate its central value to the square norm $||\cL(f,f')_{\xi}^{\chi}(g_0)||^2$.

We now recall that $\cL(f,f')$ belongs to $\cS_{\kappa}^D$~(cf.~Section 1.1). 
Before introducing the $L$-function of convolution type, we define the spinor $L$-function for a Hecke eigenform $F\in \cS_{\kappa}^D$. 
In [27,~Section 5.1] we introduced three Hecke operator $\cT_p^i$ with $0\le i\le 2$ for $p\nmid d_B$. 
Let $\Lambda_p^i$ be the Hecke eigenvalue of $\cT_p^i$ for $F$ with $0\le i\le 2$. For $p\nmid d_B$ we put
\[
Q_{F,p}(t):=1-p^{-\frac{3}{2}}\Lambda_p^1t+p^{-2}(\Lambda_p^2+p^2+1)t^2-p^{-\frac{3}{2}}\Lambda_p^1t^3+t^4.
\]
For this we note that $Q_{F,p}(p^{-s})^{-1}$ coincides with the local spinor 
$L$-function for an unramified principal series representation of the group of $\Q_p$-rational points for the split symplectic $\Q$-group $GSp(2)$ of degree two with similitudes. 
This comes from the denominator of the formal Hecke series for $GSp(1,1)_{\Q_p}\simeq GSp(2)_{\Q_p}$. 
Here recall that $GSp(2)$ is defined by
\[
GSp(2)_{\Q}:=\left\{g\in GL(4)_{\Q} \left|~{}^tg
\begin{pmatrix}
0_2 & 1_2\\
-1_2 & 0_2
\end{pmatrix}g=\nu(g)\begin{pmatrix}
0_2 & 1_2\\
-1_2 & 0_2
\end{pmatrix},~\nu(g)\in\Q^{\times}\right.\right\}.
\]
On the other hand, in [27,~Section 5.2], we introduced two Hecke operators $\cT_p^i$ with $0\le i\le 1$ for $p|d_B$. 
Let ${\Lambda'}_p^i$ be the Hecke eigenvalue of $\cT_p^i$ for $F$ with $0\le i\le 1$. For $p|d_B$ we put
\[
Q_{F,p}(t):=
\begin{cases}
(1-p^{-\frac{3}{2}}({\Lambda'}_p^1-(p-1){\Lambda'}_p^0)t+t^2)(1-{\Lambda'}_p^0p^{-\frac{1}{2}}t)&(p|\frac{d_B}{D}),\\
(1+{\Lambda'}_p^0p^{-\frac{1}{2}}t)(1-{\Lambda'}_p^0p^{-\frac{1}{2}}t)&(p|D).
\end{cases}
\]
The first one is due to Sugano [36,~(3-4)]. 
The first factor for the case $p| D$ comes from the numerator of the formal Hecke series for $GSp(1,1)_{\Q_p}$.

We define the spinor $L$-function $L(F,\operatorname{spin},s)$ of $F$ by
\[
L(F,\operatorname{spin},s):=\prod_{p\le\infty}L_p(F,\operatorname{spin},s),
\]
where
\[
L_p(F,\operatorname{spin},s):=
\begin{cases}
Q_{F,p}(p^{-s})^{-1}&(p<\infty),\\
\Gamma_{\C}(s+\frac{\kappa-1}{2})\Gamma_{\C}(s+\frac{\kappa+1}{2})&(p=\infty).
\end{cases}
\]
This is a modification of the definition in [27,~Section 5.3]. 
We can then reformulate [27,~Corollary 5.3] as follows:
\begin{prop}
The spinor $L$-function for $\cL(f,f')$ decomposes into
\[
L(\cL(f,f'),\operatorname{spin},s)=L(\pi(f),s)L(\JL(\pi(f')),s).
\]
\end{prop}
\begin{pf}
This is deduced from Lemma 2.1, Lemma 2.2 and Lemma 2.3 and the commutation relation of Hecke operators satisfied by Arakawa lifting~(cf.~[27,~Theorem 5.1]). 
Here we note that $(1-{\Lambda'}_p^0p^{-\frac{1}{2}-s})^{-1}=(1-\pi''_p(\varpi_{B,p})p^{-\frac{1}{2}-s})^{-1}=(1-\delta_p(p)p^{-\frac{1}{2}-s})^{-1}$~(respectively~$(1+{\Lambda'}_p^0p^{-\frac{1}{2}-s})^{-1}=(1+\delta_p(p)p^{-\frac{1}{2}-s})^{-1}$) holds for $p|d_B$~(respectively~$p|D$), which follows from [27,~Theorem 5.1 (ii)]. 
This is the local $L$-factor of $\JL(\pi(f'))$ at $p|d_B$~(respectively~$\pi(f)$ at $p|D$), which is a local $L$-function of a special representation~(cf.~Lemma 2.2). \qed
\end{pf}
\noindent
Of course, we see that $L(\cL(f,f'),\operatorname{spin},s)$ has the meromorphic continuation and satisfies the functional equation between $s$ and $1-s$ since so do $L(\pi(f),s)$ and $L(\JL(\pi(f')),s)$. 

We now introduce the global $L$-function
\[
L(F,\chi^{-1},s):=\prod_{p\le\infty}L_p(F,\chi^{-1},s)
\]
of convolution type for a Hecke eigenform $F\in\cS_{\kappa}^D$ and $\chi^{-1}$. 
Here the local factors $L_p(F,\chi^{-1},s)$ are given as
\[
L_p(F,\chi^{-1},s):=
\begin{cases}
Q_{F,p}(\alpha_p^{\chi}p^{-s})^{-1}Q_{F,p}(\beta_p^{\chi}p^{-s})^{-1}&(\text{$\chi$ is unramified at $p<\infty$}),\\
1&(\text{$\chi$ is ramified at $p<\infty$}),\\
\Gamma_{\C}(s+\kappa-\frac{1}{2})\Gamma_{\C}(s+\frac{1}{2})\Gamma_{\C}(s+\kappa+\frac{1}{2})\Gamma_{\C}(s+\frac{1}{2})&(p=\infty),
\end{cases}
\]
where
\[
(\alpha_p^{\chi},\beta_p^{\chi}):=
\begin{cases}
(\chi_p(\varpi_{p,1})^{-1},\chi_p(\varpi_{p,2})^{-1})&(\text{$p$:~split}),\\
(\chi_p(p)^{-1},-\chi_p(p)^{-1})=(1,-1)&(\text{$p$:~inert}),\\
(\chi_p(\varpi_p)^{-1},0)&(\text{$p$:~ramified})
\end{cases}
\]
with prime elements $\varpi_{p,1},~\varpi_{p,2}$ and $\varpi_p$ of $E_p$ introduced in Lemma 2.1 (3). 
\begin{prop}
We have
\[
L(\cL(f,f'),\chi^{-1},s)=L(\Pi,\chi^{-1},s)L(\Pi',\chi^{-1},s).
\]
This is an entire function of $s$ and satisfies the functional equation 
\[
L(\cL(f,f'),\chi^{-1},s)=\epsilon(\Pi,\chi^{-1})\epsilon(\Pi',\chi^{-1})L(\cL(f,f'),\chi^{-1},1-s),
\]
where $\epsilon(\Pi,\chi^{-1})$~(respectively~$\epsilon(\Pi',\chi^{-1})$) denotes the $\epsilon$-factor of $L(\Pi,\chi^{-1},s)$~(respectively~$L(\Pi',\chi^{-1},s)$).
\end{prop}
\begin{pf}
The first statement follows from Lemma 2.1, Lemma 2.2, Lemma 2.3 and Proposition 2.9. 
According to [15,~Corollary 19.15], $L(\Pi,\chi^{-1},s)$ and $L(\Pi',\chi^{-1},s)$ are entire functions of $s$ and satisfy 
\[
L(\Pi,\chi^{-1},s)=\epsilon(\Pi,\chi^{-1})L(\Pi,\chi,1-s),\quad 
L(\Pi',\chi^{-1},s)=\epsilon(\Pi',\chi^{-1})L(\Pi',\chi,1-s),
\]
where we note that $\Pi$ and $\Pi'$ are self-dual~(cf.~Section 2.1). 
In view of Lemma 2.1, Lemma 2.2 and Lemma 2.3, we see that $L(\Pi,\chi,s)=L(\Pi,\chi^{-1},s)$ and $L(\Pi',\chi,s)=L(\Pi',\chi^{-1},s)$, which implies the second assertion.\qed
\end{pf}
Proposition 2.10 obviously implies that the function $L(\cL(f,f'),\chi^{-1},s)$ is regular at $s=\frac{1}{2}$. 
We are now able to reformulate Theorem 2.8 as follows:
\begin{thm}
Let the assumption and the notation be as in Theorem 1.3. 
We have
\[
\frac{||\cL(f,f')_{\xi}^{\chi}(g_0)||^2}{\langle f,f\rangle\langle f',f'\rangle}=C(f,f',\xi,\chi)L(\cL(f,f'),\chi^{-1},\frac{1}{2}).
\]
\end{thm}
\subsection{Strictly positive central L-values.}
As an application of our main results we show the existence of the strictly positive central values for the $L$-functions in our concern. 
In this subsection we fix a quaternion algebra $B$ and a maximal order $\fO$ of $B$ as
\begin{align*}
&B=\Q+\Q\cdot i+\Q\cdot j+\Q\cdot k\quad(i^2=j^2=-1,~ij=-ji=k),\\
&\fO=\Z\cdot\frac{1+i+j+k}{2}+\Z\cdot i+\Z\cdot j+\Z\cdot k.
\end{align*}
We note that $d_B=2$ for this $B$. 
In [28,~Section 14] we have shown the following:
\begin{prop}
Suppose that $\xi=i/2$, which is primitive~(cf,~[28,~Section 14.2]), and $\chi$ is unramified at every finite prime, i.e. $i_p(\chi)=0$ for any finite prime $p$. 
Let $D\in\{1,2\}$ and $\kappa\ge 12$~(respectively~$\kappa\ge 8$) be divisible by $4$~(respectively~by $8$) when $D=1$~(respectively~$D=2$). 
Then there exist a primitive form $f\in S_{\kappa}(D)$ and a Hecke eigenform $f'\in\cA_{\kappa}$ such that $P_{\chi}(f,\gamma_0)P_{\chi}(f',\gamma'_0)\not=0$~(which implies 
$\cL(f,f')_{\xi}^{\chi}\not\equiv 0$).
\end{prop}
For this proposition we remark that $f$ is not assumed to be primitive in [28,~Section 14]. 
However, when $D=2$ and $\chi$ is unramified at every finite prime, we find in [28,~Section 14] $f$ satisfying $P_{\chi}(f;\gamma_0)\not=0$ and the sign condition in Section 1.3 (1). There such $f$'s are given by the powers of the primitive form of level 2 and weight 8. 
This implies that the statement remains valid even if we assume that $f$ is primitive. 
From Proposition 2.6, Proposition 2.7 and Proposition 2.12 we deduce the following:
\begin{thm}
Under the assumption in Proposition 2.12 there are Hecke eigenforms $(f,f')\in S_{\kappa}(D)\times\cA_{\kappa}$ such that
\[
L(\Pi,\chi^{-1},\frac{1}{2})>0,\quad L(\Pi',\chi^{-1},\frac{1}{2})>0.
\]
\end{thm}
\begin{pf}
Due to the isomorphism between $\cA_{\kappa}$ and the space spanned by primitive forms in $S_{\kappa+2}(d_{B})$~(cf.~Section 2.1) we have $\JL(f')\not=0$ for a non-zero $f'$. 
Furthermore we note that Proposition 2.5 implies that $L(\pi(f),\Ad,1)$ and $L(\JL(\pi(f')),\Ad,1)$, which appear in $C(f,\chi)$ and $C(f',\chi)$, are positive for a non-zero primitive form $f$ and a non-zero Hecke eigenform $f'$. 
Proposition 2.6 and Proposition 2.7 imply that $C(f,\chi)>0$ and $C(f',\chi)>0$ if $C(f,\chi)\not=0$ and $C(f',\chi)\not=0$. 
 
Let $(f,f')$ be Hecke eigenforms as in Proposition 2.12. 
In view of the explicit formulas in Proposition 2.6 and Proposition 2.7, we see $C(f,\chi)\not=0$ and $C(f',\chi)\not=0$ and thus $C(f,\chi)>0$ and $C(f',\chi)>0$ are satisfied. 
We then verify the simultaneous positivity of the two central $L$-values in the assertion.
\qed
\end{pf}
As an immediate consequence of this and Proposition 2.10 we have obtained the following theorem:
\begin{thm}
Under the same assumption as in Proposition 2.12 
there exist Hecke eigenforms $(f,f')\in S_{\kappa}(D)\times\cA_{\kappa}$ such that
\[
L(\cL(f,f'),\chi^{-1},\frac{1}{2})>0.
\]
\end{thm}
\section{Proof of Proposition 2.7.}
\subsection{}
In order to verify Theorem 2.8, which leads to Theorem 2.11, it remains to show Proposition 2.7. 
For the proof we need two formulas. 
The first one~(respectively~the second one) is due to Waldspurger~[38,~Proposition 7]~(respectively~Macdonald~[23,~Chap.V,~Section 3,~(3.4)]). 

To state the first formula we remark that every $\pi''_p$ is unitary and thus equipped with a unitary inner product. 
When $p\nmid d_B$ we embed $E_p^{\times}$ into $B_p^{\times}=GL_2(\Q_p)$ by $\iota_{\xi}$~(cf.~Section 1.4). 
Recall that, for the quadratic character $\eta$ attached to the quadratic extension $E/\Q$, we have let $L(\eta,s)$ be the $L$-function defined by $\eta$ and $L_p(\eta_p,s)$ the local factor of $L(\eta,s)$ at a place $p$~(cf.~Section 2.3). 

For the toral integral $P_{\chi}(f';b)~(b\in B_{\A_{\Q}}^{\times})$ we note that our normalized measure $ds$ of $\A_{E}^{\times}$ is the multiple of that of [38] by  $\displaystyle\frac{\sqrt{|d_{\xi}|}}{4}$. 
In fact, the volume of $\R^{\times}_+ E^{\times}\backslash\A_{E}^{\times}$ with respect to our measure~(respectively~the measure of [38]) is $\displaystyle\frac{{\bf h}(E)}{{\bf w}(E)}$~(respectively~$2L(\eta,1)=\displaystyle\frac{4{\bf h}(E)}{\sqrt{|d_{\xi}|}{\bf w}(E)}$). 
For the calculation of the volume by the measure of [38], we take the normalized measure of $\A_{\Q}^{\times}$~(cf.~Section 1.5) into account. 

For the norm $\langle f',f'\rangle$ of $f'$ we remark that our normalized measure $db$ of $B_{\A_{\Q}}^{\times}$~(cf.~Section 2.4) is the $2^{-4}3^{-1}\cdot\prod_{p|d_B}(p-1)$-multiple of that of [38], which chooses the Tamagawa measure of $B^{\times}/Z'\simeq SO(3)$. 
Actually, from the fact that the Tamagawa number of $SO(3)$ is two, we deduce that the volume of $\prod_{p<\infty}\mO_p^{\times}/\Z_p^{\times}\times\Hl^1/\{\pm1\}$ with respect to the measure of [38] is equal to $2^{3}\cdot3\cdot\prod_{p|d_B}(p-1)^{-1}$~(cf.~[20,~Section 3,~3,~(4)]), where note that $2^{-3}3^{-1}\prod_{p|d_B}(p-1)$ is the familiar number appearing in the well-known Mass formula for $B$. 
On the other hand, we see that the volume of $\prod_{p<\infty}\mO_p^{\times}/\Z_p^{\times}\times\Hl^1/\{\pm1\}$ is $\displaystyle\frac{1}{2}$ by our measure. 
These justify our remark on the measure $db$. 

Let us recall that $f'_p$ denotes the $p$-component of $f'$ for each prime $p$~(cf.~Section 2.1). With the normalized measures $ds$ and $db$ we quote [38,~Proposition 7] as follows:
\begin{prop}[Waldspurger]
For $h'=(h'_p)_{p\leq\infty}\in B^{\times}_{\A_{\Q}}$,
\[
\frac{||P_{\chi}(f';h')||^2}{\langle f',f'\rangle}=\frac{\pi\sqrt{|d_{\xi}|}}{\prod_{p|d_B}(p-1)}\frac{L(\Pi',\chi^{-1},\frac{1}{2})}{L(\pi(\JL(f')),\Ad,1)}\prod_{p\leq\infty}\alpha_p(f'_p,\chi_p,h'_p),
\]
where
\[
\alpha_p(f'_p,\chi_p,h'_p):=\frac{L_p(\eta_p,1)L_p(\pi'_p,\Ad,1)}{\zeta_p(2)L_p(\Pi'_p,\chi_p^{-1},\frac{1}{2})}\int_{\Q_p^{\times}\backslash E_p^{\times}}\frac{\langle\pi''_p(s_ph'_p)f'_p,\pi''_p(h'_p)f'_p\rangle_p}{\langle f'_p,f'_p\rangle_p}\chi_p(s_p)^{-1}ds_p
\]
with an inner product $\langle *,*\rangle_p$ of $\pi''_p$ for each place $p$.
\end{prop}
We put
\[
\phi_p(g):=\frac{\langle\pi''_p(g)f'_p,f'_p\rangle_p}{\langle f'_p,f'_p\rangle_p}\quad(g\in B_p^{\times})
\]
for each place $p$. 
We note that, at a finite prime $p$ not dividing $d_B$, $\pi''_p$ is isomorphic to $\pi'_p$ and is a unitary unramified principal series representation of $GL_2(\Q_p)$~(cf.~Section 2.2).  
The Satake parameter of $\pi''_p$ is thus of the form $(\alpha_p,\alpha_p^{-1})\in(\C^{\times})^2$ with $|\alpha_p|=1$. 
The representation $\pi''_p$ has $f'_p$ as a spherical vector. 
We remark that $\phi_p(g)$ is therefore the zonal spherical function on $GL_2(\Q_p)$ with $\phi_p(1)=1$ for $p\nmid d_B$~(for the definition see [23,~p.162]). 
The following proposition is the second formula we need~(cf.~[23,~Chap.V,~Section 3,~(3.4)]):
\begin{prop}[Macdonald]
Let $p$ be a finite prime not dividing $d_B$. 
For $a_m:=
\begin{pmatrix}
p^m & 0\\
0 & 1
\end{pmatrix}$ with a non-negative integer $m$ we have
\[
\phi_p(a_m)=\frac{p^{-\frac{m}{2}}}{1+p^{-1}}(\alpha_p^m\frac{1-p^{-1}\alpha_p^{-2}}{1-\alpha_p^{-2}}+\alpha_p^{-m}\frac{1-p^{-1}\alpha_p^2}{1-\alpha_p^2}).
\]
\end{prop}
\subsection{Calculation at $p\nmid d_B$.}
In this subsection we assume that $p\nmid d_B$. 
In what follows, we put $\mO_{E_p,i}:=\Z_p+p^i\mO_{E_p}$ for a non-negative integer $i$, and $\lambda_p:=p^{\frac{1}{2}}(\alpha_p+\alpha_p^{-1})$, which is the eigenvalue of the Hecke operator defined by the double coset $GL_2(\Z_p)
\begin{pmatrix}
p & 0\\
0 & 1
\end{pmatrix}GL_2(\Z_p)$~(cf.~[3,~Proposition 4.6.6]). 

The aim of this subsection is to prove the following proposition.
\begin{prop}
For $p\nmid d_B$,
\[
\alpha_p(f'_p,\chi_p,\gamma'_{0,p})=
\begin{cases}
1&(i_p(\chi)=0),\\
p^{-i_p(\chi)}L_p(\eta_p,1)^2&(i_p(\chi)>0),
\end{cases}
\]
where see Section 1.5 for $\gamma'_{0,p}$.
\end{prop}
\begin{pf}
\subsection*{(1)~The case of a split prime $p$.}
For the proof we do some preparation. 
For this case we note that there is an isomorphism $E_p\simeq \Q_p\times\Q_p$. 
We may thus assume that $\theta$~(cf.~Section 1.4) corresponds to $(1,0)$ via this isomorphism. 
Namely $\theta$ satisfies
\[
N(\theta)=0,~\Tr(\theta)=1,~\theta^2=\theta.
\]
We see that
\[
{\gamma'}_{0,p}^{-1}\iota_{\xi}(x+y\theta)\gamma'_{0,p}=
\begin{pmatrix}
x & 0\\
p^{-i_p(\chi)}y & x+y
\end{pmatrix}.
\]
Recall that $\varpi_{p,1}$ and $\varpi_{p,2}$ denote the two distinct prime elements of $E_p$ as in Lemma 2.1 (3). 
We may assume that $\varpi_{p,1}$ and $\varpi_{p,2}$ correspond to $(p,1)$ and $(1,p)$ respectively. Therefore
\[
\varpi_{p,1}=1+(p-1)\theta,\quad\varpi_{p,2}=1+(p-1)\bar{\theta}.
\]
For $i>0$ we have
\[
\Q_p^{\times}\backslash E_p^{\times}/\mO_{E_p,i}^{\times}\simeq\bigsqcup_{n\in\Z}\varpi_{p,1}^n(\Z_p^{\times}\backslash\mO_{E_p}^{\times}/\mO_{E_p,i}^{\times})
\]
and
\[
\Z_p^{\times}\backslash\mO_{E_p}^{\times}/\mO_{E_p,i}^{\times}\simeq\{1+b\theta\mid b\in\Z_p/p^{i}\Z_p,~1+b\in\Z_p^{\times}\}.
\]
Noting that the conductor of $\chi_p$ is $p^{i_p(\chi)}$, we can show the following lemma.
\begin{lem}
Let $p$ be split and $i_p(\chi)>0$.\\
(1)~For $k\le i_p(\chi)$,
\[
\sum_{\scriptstyle b\in p^k\Z_p/p^{i_p(\chi)}\Z_p \atop 
\scriptstyle 1+b\in\Z_p^{\times}}\chi_p(1+b\theta)^{-1}=
\begin{cases}
0&(k\le i_p(\chi)-1),\\
1 & (k=i_p(\chi)).
\end{cases}
\]
(2)~When $i_p(\chi)>1$,
\[
\sum_{\scriptstyle b\in\Z_p/p^{i_p(\chi)}\Z_p \atop 
\scriptstyle 1+b\in\Z_p^{\times},~\ord_p(b)=k}
\chi_p(1+b\theta)^{-1}=
\begin{cases}
0&(k<i_p(\chi)-1),\\
-1 & (k=i_p(\chi)-1),\\
1 & (k=i_p(\chi)).
\end{cases}
\]
When $i_p(\chi)=1$,
\[
\sum_{\scriptstyle b\in\Z_p/p\Z_p \atop 
\scriptstyle 1+b\in\Z_p^{\times},~\ord_p(b)=k}\chi_p(1+b\theta)^{-1}=
\begin{cases}
-1&(k=0),\\
1&(k=1).
\end{cases}
\]
\end{lem}
The next lemma is verified in view of the normalization of our measures of $E_p^{\times}$ and $\Z_p^{\times}$~(cf.~Section 1.5).
\begin{lem}
For a split prime $p$ we have
\[
\vol(\Z_p^{\times}\backslash\mO_{E_p,i}^{\times})=
\begin{cases}
1&(i=0),\\
p^{-i}L_p(\eta_p,1)&(i>0).
\end{cases}
\]
\end{lem}
We first carry out the proof for the case of $i_p(\chi)=0$. 
It suffices to show that
\[
\int_{\Q_p^{\times}\backslash E_p^{\times}}\frac{\langle\pi''_p(\iota_{\xi}(s_p)\gamma'_{0,p})f'_p,\pi''_p(\gamma'_{0,p})f'_p\rangle_p}{\langle f'_p,f'_p\rangle_p}\chi_p(s_p)^{-1}ds_p=\frac{\zeta_p(2)L_p(\Pi'_p\otimes\chi_p^{-1},\frac{1}{2})}{L_p(\eta_p,1)L_p(\pi'_p,\Ad,1)}.
\]
Note that we can regard $\chi_p$ as a character of $(\Q_p^{\times})^2\simeq E_p^{\times}$, thus we may write $\chi_p$ as
\[
\chi_p((t_1,t_2))=\omega_1(t_1)\omega_2(t_2)=\omega_1(t_1t_2^{-1})\quad((t_1,t_2)\in(\Q_p^{\times})^2)
\]
with two unramified characters $\omega_1$ and $\omega_2$ of $\Q_p^{\times}$ such that $\bar{\omega_1}=\omega_2$. The integral is then equal to
\[
\int_{\Q_p^{\times}\backslash E_p^{\times}}\phi_p({\gamma'}_{0,p}^{-1}\iota_{\xi}(s_p)\gamma'_{0,p})\chi_p(s_p)^{-1}ds_p=\int_{\Q_p^{\times}}\phi_p(
\begin{pmatrix}
t_1 & 0\\
0 & 1
\end{pmatrix})\omega_1(t_1)^{-1}dt_1,
\]
where we note that ${\gamma'}_{0,p}^{-1}\iota_{\xi}(s_p)\gamma'_{0,p}$ can be replaced by 
$\begin{pmatrix}
t_1 & 0\\
0 & t_2
\end{pmatrix}$ with $(t_1,t_2)\in(\Q_p^{\times})^2$ in view of the elementary divisor theorem. Now, using Proposition 3.2, we verify that this is equal to
\begin{align*}
&\frac{1}{1+p^{-1}}\sum_{m\in\Z}((p^{-\frac{1}{2}}\alpha_p)^{|m|}\frac{1-p^{-1}\alpha_p^{-2}}{1-\alpha_p^{-2}}\overline{\omega_1(p)}^m+(p^{-\frac{1}{2}}\alpha_p^{-1})^{|m|}\frac{1-p^{-1}\alpha_p^2}{1-\alpha_p^2}\overline{\omega_1(p)}^m)\\
&=\frac{1}{1+p^{-1}}\{(-1+(1-\alpha_pp^{-\frac{1}{2}}\overline{\omega_1(p)})^{-1}+(1-\alpha_pp^{-\frac{1}{2}}\overline{\omega_2(p)})^{-1})\frac{1-p^{-1}\alpha_p^{-2}}{1-\alpha_p^{-2}}\\
&+(-1+(1-\alpha_p^{-1}p^{-\frac{1}{2}}\overline{\omega_1(p)})^{-1}+(1-\alpha_p^{-1}p^{-\frac{1}{2}}\overline{\omega_2(p)})^{-1})\frac{1-p^{-1}\alpha_p^{2}}{1-\alpha_p^{2}}\}.
\end{align*}
By a direct calculation we prove that this coincides with the ratio of the local L-functions on the right hand side~(cf.~Lemma 2.1). 
We thus have $\alpha_p(f'_p,\chi_p,\gamma'_{0,p})=1$.

We next consider the case of $i_p(\chi)>0$. 
Since $\phi_p$ is bi-invariant by $GL_2(\Z_p)$ and trivial on the center we see that, for $1+b\theta\in O_{E_p}^{\times}$ and $n\in\Z\setminus\{0\}$, 
\begin{align*}
\phi_p({\gamma'}_{0,p}^{-1}\iota_{\xi}(1+b\theta)\gamma'_{0,p})&=
\phi_p(
\begin{pmatrix}
p^{2(i_p(\chi)-\ord_p(b))} & 0\\
0 & 1
\end{pmatrix}),\\
\phi_p({\gamma'}_{0,p}^{-1}\iota_{\xi}(\varpi_1^n(1+b\theta))\gamma'_{0,p})&=\phi_p(
\begin{pmatrix}
p^{2i_p(\chi)+|n|} & 0\\
0 & 1
\end{pmatrix}).
\end{align*}
By Proposition 3.2, Lemma 3.4 and Lemma 3.5 the local integral $\displaystyle\int_{\Q_p^{\times}\backslash E_p^{\times}}\phi_p({\gamma'}_{0,p}^{-1}\iota_{\xi}(s_p)\gamma'_{0,p})\chi_p(s_p)^{-1}ds_p$ involved in $\alpha_p(f'_p,\chi_p,\gamma'_{0,p})$ is
\begin{align*}
&\vol(\mO_{E_p,i_p(\chi)}^{\times}/\Z_p^{\times})\{\sum_{n\in\Z\setminus\{0\}}\phi_p(
\begin{pmatrix}
p^{2i_p(\chi)+|n|} & 0\\
0 & 1
\end{pmatrix})\chi(\varpi_1^n)^{-1}\sum_{\scriptstyle b\in Z_p/p^{i_p(\chi)}\Z_p \atop \scriptstyle 1+b\in\Z_p^{\times}}\chi(1+b\theta)^{-1}\\
&+\sum_{\scriptstyle b\in\Z_p/p^{i_p(\chi)}\Z_p \atop \scriptstyle 1+b\in\Z_p^{\times}}\phi_p(
\begin{pmatrix}
p^{2(i_p(\chi)-\ord_p(b))} & 0\\
0 & 1
\end{pmatrix})\chi(1+b\theta)^{-1}\}\\
=&p^{-i_p(\chi)}L_p(\eta_p,1)(1-\phi_p(
\begin{pmatrix}
p^2 & 0\\
0 & 1
\end{pmatrix}))=\frac{L_p(\eta_p,1)((p+1)^2-\lambda_p^2)}{p^{i_p(\chi)}\cdot p(p+1)}.
\end{align*}
On the other hand, Lemma 2.1 implies that the ratio of the local $L$-functions in $\alpha_p(f'_p,\chi_p,\gamma'_{p,0})$ is equal to
\[
\frac{L_p(\eta_p,1)\cdot p(p+1)}{(p+1)^2-\lambda_p^2}.
\]
As a result we deduce that
\[
\alpha_p(f'_p,\chi_p,\gamma'_{p,0})=p^{-i_p(\chi)}L_p(\eta_p,1)^2.
\]
\subsection*{(2)~The case of an inert prime $p$.}
A complete set of representatives for $\Q_p^{\times}\backslash E_p^{\times}/\mO_{E_p,i}^{\times}\simeq\Z_p^{\times}\backslash \mO_{E_p}^{\times}/\mO_{E_p,i}^{\times}$ is given by
\[
\{1+b\theta\mid b\in\Z_p/p^{i}\Z_p\}\sqcup\{a+\theta\mid a\in p\Z_p/p^{i}\Z_p\}
\]
for $i>0$. 
To verify the proposition for this case we need the following two lemmas, whose  proofs are similar to those of Lemma 3.4 and Lemma 3.5.
\begin{lem}
Let $p$ be inert and $i_p(\chi)>0$.\\
(1)~We have
\[
\sum_{b\in\Z_p/p^{i_p(\chi)}\Z_p}\chi_p(1+b\theta)^{-1}=-\sum_{a\in p\Z_p/p^{i_p(\chi)}\Z_p}\chi_p(a+\theta)^{-1}=
\begin{cases}
-\chi(\theta)^{-1}&(i_p(\chi)=1),\\
0&(i_p(\chi)>1).
\end{cases}
\]
(2)~When $i_p(\chi)>1$,
\[
\sum_{\scriptstyle b\in\Z_p/p^{i_p(\chi)}\Z_p \atop \scriptstyle \ord_p(b)=k}\chi_p(1+b\theta)^{-1}=
\begin{cases}
0&(k<i_p(\chi)-1),\\
-1&(k=i_p(\chi)-1),\\
1&(k=i_p(\chi)).
\end{cases}
\]
\end{lem}
\begin{lem}
For an inert prime $p$ we have
\[
\vol(\Z_p^{\times}\backslash\mO_{E_p,i}^{\times})=
\begin{cases}
1&(i=0),\\
p^{-i}L_p(\eta_p,1)&(i>0).
\end{cases}
\]
\end{lem}
The local integral involved in $\alpha_p(f'_p,\chi_p,\gamma'_{0,p})$ is expressed as
\[
\vol(\Z_p\backslash\mO_{E_p,i_p(\chi)}^{\times})\int_{\Z_p^{\times}\backslash\mO_{E_p}^{\times}/\mO_{E_p,i_p(\chi)}^{\times}}\phi_p(\gamma_{0,p}^{-1}\iota_{\xi}(s_p)\gamma_{0,p})\chi_p(s_p)^{-1}ds_p.
\]
First let $i_p(\chi)=0$. 
This integral is evaluated to be $\vol(\Z_p^{\times}\backslash \mO_{E_p^{\times}})=1$. 
Let us calculate the ratio of the local L-functions occurring in $\alpha_p(f'_p,\chi_p,\gamma'_{0,p})$. 
With the Satake parameter $(\alpha_p,\alpha_p^{-1})$ of $\pi'_p$ we have
\begin{align*}
&\frac{L_p(\eta_p,1)L_p(\pi'_p,\Ad,1)}{\zeta_p(2)L_p(\Pi'_p,\chi_p^{-1},\frac{1}{2})}=\frac{(1+p^{-1})^{-1}(1-p^{-1})^{-1}(1-\alpha_p^2p^{-1})^{-1}(1-\alpha_p^{-2}p^{-1})^{-1}}{(1-p^{-2})^{-1}(1-\alpha_pp^{-\frac{1}{2}})^{-1}(1-\alpha_p^{-1}p^{-\frac{1}{2}})^{-1}(1+\alpha_pp^{-\frac{1}{2}})^{-1}(1+\alpha_p^{-1}p^{-\frac{1}{2}})^{-1}}\\
&=1
\end{align*}
(cf.~Lemma 2.1), 
which implies the assertion for this case. 

Next let $i_p(\chi)>0$. 
In view of the bi-invariance by $GL_2(\Z_p)$ and the triviality on the center for $\phi_p$ we see that 
\[
\phi_p({\gamma'_{0,p}}^{-1}\iota_{\xi}(1+b\theta)\gamma'_{0,p})=\phi_p(
\begin{pmatrix}
p^{2(i_p(\chi)-\ord_p(b))} & 0\\
0 & 1
\end{pmatrix}),\quad \phi_p({\gamma'_{0,p}}^{-1}\iota_{\xi}(a+\theta)\gamma'_{0,p})=\phi_p(
\begin{pmatrix}
p^{2i_p(\chi)} & 0\\
0 & 1
\end{pmatrix})
\]
for $a\in p\Z_p/p^{i_p(\chi)}\Z_p$ and $b\in\Z_p/p^{i_p(\chi)}\Z_p$. 
By Proposition 3.2, Lemma 3.6 and Lemma 3.7 we therefore evaluate the local integral to be $\displaystyle\frac{L_p(\eta_p,1)^2((p+1)^2-\lambda_p^2)}{p^{i_p(\chi)}\cdot p^2}$. 
On the other hand, recalling Lemma 2.1, we see that the ratio of the local L-functions in $\alpha_p(f'_p,\chi_p,\gamma'_{p,0})$ is $p^{2}((p+1)^2-\lambda_p^2)^{-1}$. 
We then obtain the evaluation of $\alpha_p(f'_p,\chi_p,\gamma'_{0,p})$ in the assertion.
\subsection*{(3)~The case of a ramified prime $p$.}
We give a complete set of representatives for $\Q_p^{\times}\backslash E_p^{\times}/\mO_{E_p,i}^{\times}$ as follows:
\[
\{1+b\theta\mid b\in\Z_p/p^{i}\Z_p\}\sqcup\{ap+\theta\mid a\in \Z_p/p^{i}\Z_p\},
\]
where $i\ge0$. 
We state the following two lemmas, whose proofs are similar to those of the corresponding two lemmas for the two cases above.
\begin{lem}
Let $p$ be ramified and $i_p(\chi)>0$.\\
(1)~We have
\[
\sum_{b\in\Z_p/p^{i_p(\chi)}\Z_p}\chi_p(1+b\theta)^{-1}=-\sum_{a\in \Z_p/p^{i_p(\chi)}\Z_p}\chi_p(ap+\theta)^{-1}=0.
\]
(2)~When $i_p(\chi)>1$,
\[
\sum_{\scriptstyle b\in\Z_p/p^{i_p(\chi)}\Z_p \atop \scriptstyle \ord_p(b)=k}\chi_p(1+b\theta)^{-1}=
\begin{cases}
0&(k<i_p(\chi)-1),\\
-1&(k=i_p(\chi)-1),\\
1&(k=i_p(\chi)).
\end{cases}
\]
\end{lem}
\begin{lem}
For a ramified prime $p$ we have
\[
\vol(\Z_p^{\times}\backslash\mO_{E_p,i}^{\times})=
\begin{cases}
1&(i=0),\\
p^{-i}=p^{-i}L_p(\eta_p,1)&(i>0).
\end{cases}
\]
\end{lem}
Suppose first that $i_p(\chi)=0$. The local integral is
\[
\int_{\Z_p^{\times}\backslash\mO_{E_p}^{\times}}\phi_p({\gamma'}_{0,p}^{-1}\iota_{\xi}(s_p)\gamma'_{0,p})\chi_p(s_p)^{-1}ds_p+\int_{\Z_p^{\times}\backslash\varpi_{p}\mO_{E_p}^{\times}}\phi_p({\gamma'}_{0,p}^{-1}\iota_{\xi}(s_p)\gamma'_{0,p})\chi_p(s_p)^{-1}ds_p,
\]
where recall that $\varpi_p$ denotes a prime element of $E_p$~(cf.~Lemma 2.1 (3)). 
The first integral is proved to be $1$. 
By virtue of the elementary divisor theorem we see that the second integral is reduced to $\phi_p(
\begin{pmatrix}
p & 0\\
0 & 1
\end{pmatrix})\chi_p(\varpi_{p})^{-1}$, which is equal to
\[
\frac{p^{-\frac{1}{2}}}{1+p^{-1}}(\alpha_p\frac{1-p^{-1}\alpha_p^{-2}}{1-\alpha_p^{-2}}+\alpha_p^{-1}\frac{1-p^{-1}\alpha_p^2}{1-\alpha_p^2})\chi_p(\varpi_{p})^{-1}=\frac{\lambda_p\omega_p(p)}{p+1}.
\]
Here see Section 2.2 for $\omega_p$. Thus we have
\[
\int_{\Q_p^{\times}\backslash E_p^{\times}}\phi_p(\gamma_{0,p}^{-1}\iota_{\xi}(s_p)\gamma_{0,p})\chi_p(s_p)^{-1}ds_p=\frac{p+1+\lambda_p\omega_p(p)}{p+1}.
\]
On the other hand, the ratio of the local L-functions in $\alpha_p(f'_p,\chi_p,\gamma'_{0,p})$ is verified to be
\[
\frac{(p+1)(p+1-\lambda_p\omega_p(p))}{(p+1)^2-\lambda_p^2}.
\]
The assertion is now obvious.

Next let $i_p(\chi)>0$. 
We have 
\[
\phi_p({\gamma'}_{0,p}^{-1}\iota_{\xi}(1+b\theta)\gamma'_{0,p})=\phi_p(
\begin{pmatrix}
p^{2(i_p(\chi)-\ord_p(b))} & 0\\
0 & 1
\end{pmatrix}),\quad \phi_p({\gamma'}_{0,p}^{-1}\iota_{\xi}(ap+\theta)\gamma'_{0,p})=\phi_p(
\begin{pmatrix}
p^{2i_p(\chi)+1} & 0\\
0 & 1
\end{pmatrix})
\]
for $a\in \Z_p/p^{i_p(\chi)}\Z_p$ and $b\in\Z_p/p^{i_p(\chi)}\Z_p$. 
By Proposition 3.2, Lemma 3.8 and Lemma 3.9 we evaluate the local integral to be 
$\displaystyle\frac{(p+1)^2-\lambda_p^2}{p^{i_p(\chi)}\cdot p(p+1)}$. 
The ratio of the local $L$-functions is $\displaystyle\frac{L_p(\pi'_p,\Ad,1)}{\zeta_p(2)}=\displaystyle\frac{p(p+1)}{(p+1)^2-\lambda_p^2}$. 
The formula for this case then follows immediately. 
As a result we have completed the proof of the proposition.\qed
\end{pf}
\subsection{Calculation at $p|d_B$.}
Throughout this subsection we assume that $p|d_B$. Then $p$ is never split in $E/\Q$. 
Recall that we have assumed that $\chi_p$ is unramified at such $p$~(cf.~Assumption 1.2). 
We now show the following:
\begin{prop}
Recall that $r_p$ denotes the ramification index of $p$ for the quadratic extension $E/\Q$~(cf.~Section 2.4). 
We have
\[
\alpha_p(f'_p,\chi_p,\gamma'_{0,p})=
\begin{cases}
r_p(1-p^{-1})&(\text{$p$ is inert, or $p$ is ramified and $\pi''_p|_{E_p^{\times}}=\chi_p$}),\\
0&(\text{$p$ is ramified and $\pi''_p|_{E_p^{\times}}\not=\chi_p$}),
\end{cases}
\]
where we note that $\gamma'_{0,p}=\varpi_{B,p}$~(cf.~Section 1.5). 
\end{prop}
Our proof of this needs the following:
\begin{lem}
(1)
\[
L_p(\pi'_p,\Ad,1)=(1-p^{-2})^{-1}.
\]
(2)
\[
L_p(\Pi'_p,\chi_p^{-1},\frac{1}{2})=
\begin{cases}
(1-p^{-2})^{-1}&(\text{$p$:inert})\\
(1-\delta_p(p)\omega_p(p)p^{-1})^{-1}&(\text{$p$:ramified})
\end{cases},
\]
where see Section 2.2 for $\delta_p$ and $\omega_p$.
\end{lem}
\begin{pf}
This is a direct consequence of Lemma 2.2.\qed
\end{pf}
Now we are able to calculate $\alpha_p(f'_p,\chi_p,\gamma'_{0,p})$. We first note that we can replace $\gamma'_{0,p}$ by $1$ since $\pi''_p(\varpi_{B,p})\in\{\pm1\}$. By a direct calculation we verify
\[
\int_{\Q_p^{\times}\backslash E_p^{\times}}\frac{\langle\pi''_p(s_p\gamma'_{0,p})f'_p,\pi''_p(\gamma'_{0,p})f'_p\rangle_p}{\langle f'_p,f'_p\rangle_p}\chi_p(s_p)^{-1}ds_p=
\begin{cases}
r_p&(\text{$p$ is inert, or $p$ is ramified and $\pi''_p|_{E_p^{\times}}=\chi_p$}),\\
0&(\text{$p$ is ramified and $\pi''_p|_{E_p^{\times}}\not=\chi_p$}).
\end{cases}
\]
Now note that $\pi''_p|_{E_p^{\times}}=\chi_p$ means $\delta_p=\omega_p$. 
As for the ratio of the local L-functions, Lemma 3.11 yields
\[
\frac{L_p(\eta_p,1)L_p(\pi'_p,\Ad,1)}{\zeta_p(2)L_p(\Pi'_p,\chi_p^{-1},\frac{1}{2})}=1-p^{-1}
\]
unless $p$ is ramified and satisfies $\pi''_p|_{E_p^{\times}}\not=\chi_p$. These imply Proposition 3.10.
\subsection{Calculation at $\infty$.}
To complete the proof of Proposition 2.7 we are left with the calculation of $\alpha_{\infty}(f'_{\infty},\chi_{\infty}^{-1},\gamma'_{0,\infty})$, where $\gamma'_{0,\infty}=1$. 
We note that the inner product $\langle *,*\rangle_{\infty}$ is taken as $(*,*)_{\kappa}$~(cf.~Section 2.4).
\begin{prop}
\[
\alpha_{\infty}(f'_{\infty},\chi_{\infty}^{-1},\gamma'_{0,\infty})=\frac{\kappa+1}{4\pi}\frac{(f'_{\infty,\kappa},f'_{\infty,\kappa})_{\kappa}}{(f'_{\infty},f'_{\infty})_{\kappa}},
\]
where see Section 2.4 for $f'_{\infty,\kappa}$.
\end{prop}
This proposition follows from two lemmas. 
The first one below is settled by a direct calculation:
\begin{lem}
\[
\int_{E_{\infty}^{\times}/\Q_{\infty}^{\times}}\phi_{\infty}(s_{\infty})\chi_{\infty}(s_{\infty})^{-1}ds_{\infty}=\frac{(f'_{\infty,\kappa},f'_{\infty,\kappa})_{\kappa}}{2(f'_{\infty},f'_{\infty})_{\kappa}}.
\]
\end{lem}
To state the second lemma we note that the archimedean component $\pi'_{\infty}$ of $\pi'$ is the discrete series of weight $\kappa+2$~(cf.~[35,~Section 6]). We have the following:
\begin{lem}
\begin{align*}
&L_{\infty}(\pi'_{\infty},\Ad,1)=2^{-(\kappa+1)}\pi^{-(\kappa+3)}(\kappa+1)!,\\
&L_{\infty}(\Pi'_{\infty},\chi_{\infty}^{-1},\frac{1}{2})=2^{-\kappa}\pi^{-(\kappa+2)}\kappa!,\\
&\zeta_{\infty}(2)=L_{\infty}(\eta_{\infty},1)=\pi^{-1}.
\end{align*}
\end{lem}
\noindent
The first two follow from Lemma 2.3 and the last one is well-known.

As a result of Proposition 3.1, Proposition 3.3, Proposition 3.10 and Proposition 3.12, we have proved Proposition 2.7.
\subsection*{Acknowledgement.}
We are very grateful to Kimball Martin for his comments on the explicit formulas for the toral integrals in terms of the central $L$-values for $GL(2)$. 
Our deep gratitude is due to Ralf Schmidt for a fruitful discussion.

\end{document}